\documentclass[11pt]{amsart} 

\usepackage[vmargin=3cm, hmargin=3cm]{geometry}
\usepackage[foot]{amsaddr}

\usepackage{mathpazo}

\usepackage{multirow}
\usepackage{longtable}

\usepackage[mathcal]{eucal}

\DeclareFontFamily{OT1}{pzc}{}
\DeclareFontShape{OT1}{pzc}{m}{it}{<-> s * [1.10] pzcmi7t}{}
\DeclareMathAlphabet{\mathpzc}{OT1}{pzc}{m}{it}

\usepackage{afterpage}
\usepackage{makecell}
\usepackage{changepage}
\usepackage{wrapfig}
\usepackage{enumitem}
\usepackage{hyperref}
\usepackage{mathrsfs}  
\usepackage{latexsym}
\usepackage{enumitem}
\usepackage{bbm} 
\usepackage{tikz}

\usepackage{yhmath}
\usepackage{fancyhdr}
\usepackage{amsmath, amsfonts, amssymb, amsthm, tensor}
\usepackage[all]{xy}
\usepackage{verbatim}
\usepackage{graphicx}
\usepackage{tikz-cd}
\usetikzlibrary{decorations.pathreplacing,calligraphy}
\usepackage[T1]{fontenc}
\usepackage[utf8]{inputenc}
\usepackage{multicol}
\usepackage{hhline}
\usepackage{color}
\usepackage{adjustbox}
\usepackage{forest}
\usetikzlibrary{decorations.pathreplacing,calc}
\usetikzlibrary{arrows.meta,bending,quotes}
\usetikzlibrary{arrows}

\usepackage{mathtools}

\def\co{\colon\thinspace}

\mathchardef\mhyphen="2D

\numberwithin{equation}{section}

\newtheorem{thmA}{Theorem}

\newtheorem{corA}[thmA]{Corollary}

\newtheorem{teorema}{Theorem}[section]  
\newtheorem{lem}[teorema]{Lemma} 
\newtheorem{prop}[teorema]{Proposition}
\newtheorem{cor}[teorema]{Corollary}

\theoremstyle{definition}
\newtheorem{defin}[teorema]{Definition}

\theoremstyle{remark} 
\newtheorem{ex}[teorema]{Example}
\newtheorem{oss}[teorema]{Remark}

\newcommand{\C}{\mathbb{C}}
\newcommand{\R}{\mathbb{R}}
\newcommand{\N}{\mathbb{N}}
\newcommand{\Z}{\mathbb{Z}}
\newcommand{\Q}{\mathbb{Q}}
\newcommand{\T}{\mathbb{T}}

\newcommand{\trop}{\textup{trop}}

\newcommand{\DTrop}{\mathpzc{Trop}^d}

\newcommand{\B}{\mathbb{B}}

\newcommand{\bigmid}{\, \bigg | \, }

\newcommand{\posR}{(\R_{\ge 0})^m}
\newcommand{\pkcone}[1]{K[\![#1]\!]}
\newcommand{\pk}{K[\![t]\!]}

\newcommand{\pt}{\T[\![t]\!]}
\newcommand{\pb}{\B[\![t]\!]}

\newcommand{\ptt}[1]{\T[\![t_1, \dots , t_{#1}]\!]}
\newcommand{\pbb}[1]{\B[\![t_1, \dots , t_{#1}]\!]}
\newcommand{\pkk}[1]{K[\![t_1, \dots , t_{#1}]\!]}

\newcommand{\fatg}{\textbf{g}}

\newcommand{\hh}{\textbf{h}}
\newcommand{\jj}{\textbf{j}}
\newcommand{\rr}{\textbf{r}}

\newcommand{\triv}{\textup{triv}}

\newcommand{\diff}[2]{#1\{x_1, \ldots , x_{#2}\}}
\newcommand{\basic}[2]{#1\{x_1, \ldots , x_{#2}\}_\mathrm{basic}}

\newcommand{\norm}[1]{|\!|{#1}|\!|}
\newcommand{\ndiv}{\hspace{-4pt}\not|\hspace{2pt}}
\newcommand{\initial}{\mathrm{in}}
\newcommand{\Initial}{\mathrm{In}}

\newcommand{\supp}{\mathrm{supp}}

\DeclareMathOperator*{\colim}{colim} 

\newcommand{\alg}[1]{#1\textbf{\textup{-Alg}}}

\usepackage{todonotes}

\title[Fundamental theorem of trop. diff. algebra over nontrivially valued fields]{The Fundamental theorem of tropical differential algebra over nontrivially valued fields and the radius of convergence of nonarchimedean differential equations}

\author{Francesco Gallinaro}
\email{francesco.gallinaro@sns.it}
\address{Centro di Ricerca Matematica ``Ennio De Giorgi", Scuola Normale Superiore, Piazza dei Cavalieri 3, 56125 Pisa, Italy}

\author{Stefano Mereta}
\email{stefano.mereta@cunef.edu}
\address{Departamento de Matem\'aticas, CUNEF Universidad, Calle de Almansa 101, 28040 Madrid, Spain}

\date{\today}


\begin{document}
	\begin{abstract}
        We prove a fundamental theorem for tropical
        partial differential equations, analogous to the fundamental theorem of tropical geometry in this context. We extend results from Aroca et al., Falkensteiner et al. and from Fink and Toghani for the case of trivial valuation as introduced by Grigoriev to differential equations with power series coefficients over any valued field. Crucial ingredients are the framework for tropical partial differential equations introduced by Giansiracusa and Mereta and a result on infinite intersections of projections of fibers of tropicalizations, which we prove using Hrushovski and Loeser's model-theoretic interpretation of Berkovich analytification. As a corollary of the fundamental theorem, we show that the radius of convergence of solutions of an ordinary differential equation over a nontrivially valued field can be computed tropically.
	\end{abstract}

	\maketitle
 
\setcounter{tocdepth}{1}	
\tableofcontents

\section*{Introduction}

 Differential equations are ubiquitous in the sciences from theoretical physics to ecology. While a great number of analytical and numerical methods to compute their solutions have been developed since Newton, purely algebraic methods are more recent and generally less developed. The last decade has seen a resurgence of interest in this topic. The theory is based mostly on the pioneering work of Ritt \cite{ritt} and Kolchin \cite{kolchin} on the development of differential algebra, that allows one to study systems of algebraic ordinary (and partial) differential equations by means of algebraic tools, such as Gröbner bases, and the works of Macaulay and later Gröbner \cite{grobmac} on the use of differential operators in commutative algebra. This led to the implementation of software routines computing the symbolic solutions of system of linear PDEs (see for example \cite{ridamarc}, \cite{YaironMac}).
 
Algebraic (partial) differential equations are formed from polynomial expressions in an
indeterminate function $f$ and its (partial) derivatives. This smaller class of equations comprises most of the differential equations appearing in applications and in pure mathematics. 

Tropical geometry is a fairly new area of algebraic geometry that started developing from the early 2000s and allows to move problems of algebro-geometric nature to a piecewise linear setting  via a valuation. Here methods from combinatorics and polyhedral geometry can be applied to study properties of the initial object. A comprehensive reference for tropical geometry is \cite{macsturm}.

Tropical methods for algebraic differential equations were introduced for the first time in \cite{grig}. There, the author defined tropical differential equations, their solutions, and a framework for tropicalizing algebraic ODEs over a ring of formal power series $\pk$, and their solutions. The field $K$ is endowed with the trivial valuation. In this context, one tropicalizes differential equations by recording the $t$-adic valuation
of each coefficient, and the tropicalization of power series solutions is performed by applying the trivial valuation of $K$ coefficientwise, i.e.\ recording their support. The same can also be adapted to the case of PDEs. 
Note that, when $K$ can be endowed with a nontrivial nonarchimedean valuation, Grigoriev's tropicalization does not record any information about the valuation of the coefficients in a differential equation or a power series solution.  Thus all the information about convergence of power series solutions is lost.

In order to remedy this flaw and preserve finer information about the coefficients, in \cite{gianmereta} the authors abstracted the main features of Grigoriev's setting and introduced a framework that extends and refines it. This allows one to perform a tropicalization that, in case the field $K$ has a non-trivial valuation (for example $K$ being the field of Puiseux series or that of $p$-adic numbers), preserves the valuation of the coefficients of the solutions and of the equations. While in the aformentioned paper only results for ODEs are presented, a general framework for PDEs can also be developed (see \cite{tesimereta}). It is in this language that we will phrase our results, always making clear how these look in the case of ODEs.

Analogously to the fundamental theorem of tropical geometry \cite[Theorem 3.2.3]{macsturm} (also called Kapranov's theorem in the case of hypersurfaces. See \cite{kapranov}), one can ask if solutions to a differential equation tropicalize to solutions of its tropicalization, and if all the tropical solutions can be lifted back. For Grigoriev's setting this question has  been positively answered in the case of ODEs by Aroca et
al.~in \cite{aroca} and subsequently generalised to PDEs in \cite{sebastian} (assuming $K$ is an uncountable algebraically closed field of
characteristic 0). Furthermore, paralleling the role of Gr{\"o}bner theory in the non-differential setting, where a tropical variety is the set of weight vectors whose initial ideal does not contain any monomial, the authors of \cite{FT20} and \cite{hugao} define initial forms and develop a Gr{\"o}bner-theoretic approach to Grigoriev's framework for ODEs. 

In tropical geometry the necessary datum to perform the tropicalization of a set of polynomial equations over a ring $R$ and of its solutions is that of a nonarchimedean valuation $v:R \rightarrow S$ to a semiring $S$. To introduce a similar notion of tropicalization in the differential case we will treat here, we need a \emph{differential enhancement} of such a valuation. See Definition \ref{def:diff-enhan}. For now, let us just say that if $R$ is a differential ring and $v: R \rightarrow S$ is a valuation then a differential enhancement of $v$ is a pair $\mathbf{v}=(v,\widetilde{v})$ where $\widetilde{v}$ is a map from $R$ to a space endowed with a differential structure that is somehow compatible with $v$. 

The main example is roughly the following one (see Section \ref{section:tropdiffeq} for the details). Let $R^K_m$ denote the power series ring $\pkk m$, where $K$ is a valued field with value group $\Gamma_K$. Refining \cite{sebastian}, we will introduce a valuation $v$ on $R^K_m$ which takes values in a semiring $\textup{Conv}_m(\widetilde \Gamma_K)$, where $\widetilde \Gamma_K$ is the semifield $\Gamma_K \cup \{\infty\}$, whose elements are convex hulls  of certain subsets of $\N^m$ with vertices weighted in $\Gamma_K$. We obtain a differential enhancement of $\mathbf v =(v, \widetilde v)$ by taking $\widetilde v$ to be a map from $R_m$ into a ring of formal power series with coefficients in the semifield $\widetilde \Gamma_K$, so that the diagram 
\begin{center}
\begin{tikzcd}
& \widetilde \Gamma_K  [\![t_1 , \dots , t_m]\!]  \arrow[d]  \\
R^K_m \arrow[ur, "\widetilde{v}"] \arrow[r, swap, "v"] &  \textup{Conv}_m(\widetilde \Gamma_K)
\end{tikzcd}
\end{center}
commutes, where the vertical arrow is defined by taking the weighted convex hull of the support of a series. We will tropicalize equations by means of $v$ and solutions via $\widetilde v$. The vertical arrow in the diagram above is a tropical pair, as introduced in \cite{gianmereta}, that we will denote as $\mathbf{\Gamma}^K_m$.

One of the contributions of the present work is the generalization of the theory developed in \cite{FT20, hugao} to the nontrivially valued case and to PDEs. In this case, the initial form of a differential polynomial with coefficients in $R^K_m $ will depend on a weight vector of tropical power series and on a vector $\omega \in \R^m$ with rationally independent coordinates, and the reason is that in the process of defining the initial form we will make use of the field of multivariate Laurent series with respect to such a vector $\omega$ (see \cite{arocacanojung, arocarond}).

The main goal of this paper is to prove the following theorem, which is the generalization to the nontrivially valued case of the fundamental theorems of \cite{aroca}, \cite{sebastian} and \cite{FT20}:
\begin{thmA} [Fundamental theorem of tropical differential algebra]\label{theorem:fundamental-introduction}
        Let $K$ be an algebraically closed field of characteristic 0 and $v_K : K \rightarrow \T$ a valuation. There is an algebraically closed valued field $(L,v_L)$ extending $(K,v_K)$, with a differentially enhanced valuation \[
            \mathbf{v}=(v,\widetilde{v})\co (R_m^L , \{d/dt_i\}_{i \in [m]}) \to \mathbf{\Gamma}_m^L
        \] as described above, such that for every differential ideal $I$ in $R_m^K\{x_1,\dots,x_n\}$, the following subsets of $(\widetilde \Gamma_L [\![t_1 , \dots , t_m]\!])^n$ coincide:
        \begin{enumerate}
        \item the tropicalization with respect to $\widetilde v$ of the set of solutions to $I$ in $R_m^L$;
        \item the solutions to the tropicalization with respect to $v$ of the differential ideal $IR_m^L\{x_1,\dots,x_n\}$;
        \item the set of tropical power series vectors $S$ such that the initial ideal $\Initial_{(S, \omega)}(I)$ of $I$ with respect to $(S,\omega)$ does not contain monomials for any vector $\omega \in (\R_{>0})^m$ with rationally independent coordinates.
        \end{enumerate}
	\end{thmA}
    
All the objects appearing in the statement above will be introduced in the following. Theorem \ref{theorem:fundamental-introduction} will be proved in two parts, in Section \ref{section:theorem} and in Section \ref{section:grobner}. The proof of the result above makes use of the fundamental theorem of tropical geometry and is a conceptual improvement of the proof of the fundamental theorem presented in \cite{aroca}, \cite{sebastian} and \cite{FT20}. As a step of the proof we prove an approximation statement (Proposition \ref{Prop:limit-empty-iff-some-empty}) for infinite intersections of projections of fibers of tropicalization, using methods from model theory, specifically from the Hrushovski--Loeser interpretation of Berkovich spaces as type spaces \cite{HrushovskiLoeser}. A similar approximation result is used in \cite{aroca} in the case of a trivial valuation. A major difference is that our statement finds a point in a certain infinite intersections of sets only after passing to an extension of the ambient valued field: without this extension the statement, and indeed Theorem \ref{theorem:fundamental-introduction}, become false in full generality, even assuming that the valued field is uncountable, as we show in Example \ref{ex:Puiseux}. We show that the extension is not necessary for ideals with simple structure (Proposition \ref{Prop:limit-empty-iff-some-empty} (2)). 

Finally, we introduce a notion of radius of convergence for tropical solutions in the case of ODEs. We recall that the convergence of a power series with coefficients in a non-archimedean valued field is equivalent to the limit of the norms of the coefficients being 0. As the coefficients of a tropical solution are tropical numbers, we can define a notion of convergence for such a tropical power series. Since Theorem \ref{theorem:fundamental-introduction} above states the equality between the tropicalized solutions and the tropical ones, as an important application it implies the following corollary:
\begin{corA}\label{corolla:intro}
For ODEs, if there is a tropical solution with radius of convergence $r$ then there exist a valued field extension $(L,v_L)$ of $(K,v_K)$ and a classical solution with coefficients in $L$ with radius of convergence $r$.
\end{corA}

Hence, by computing the tropical solutions we can study the radius of convergence of the classical ones, without having any knowledge about them. This paper takes the first steps towards the introduction of tropical techniques for the computation of the radius of convergence function of the solutions to ordinary differential equations over the $p$-adics (and in general over a nonarchimedean field), one of the central problems in the theory of $p$-adic differential equations. Indeed, an interesting feature of $p$-adic ODEs is that the radius of convergence of their solutions is not controlled by some "visible" object, such as, in the complex case, the poles of the coefficients of the equation: in the $p$-adic world even equations as easy as that of the exponential give solutions with finite radius of convergence at any point. The topology of the space itself is an obstacle to the convergence. 
The language of Berkovich geometry, introduced in \cite{ber}, has proved to be the right one to describe phenomena related to these radii of convergence. The radius of convergence of the solutions of a linear $p$-adic ODE on a Berkovich curve as a function of the expansion point has been proved to be a continuous piecewise linear function (see \cite{ christol1994modules, baldassarri2010continuity, baldassarri2007continuity, kedlaya}) with finitely many changes of slope \cite{christol2011radius, pulita2015convergence, pulita2015continuity,pulita2015convergence2, poineau2013convergence}, whose behaviour is actually controlled by a finite skeleton on which the curve retracts on. In general, although an explicit iterative formula to compute the radius of converge exists (see \cite{christol2011radius}), it is difficult to calculate it. A comprehensive reference for the subject is \cite{kedlaya}. See also \cite{librodwork, librorobba}.

\subsection*{Organisation of the paper} In Section \ref{section:tropdiffeq} we recall definitions and results of the theory of tropical differential equations, mainly in the language of \cite{gianmereta}, encompassing the theory of Grigoriev. In Section \ref{section:preliminaries} we prove Proposition \ref{Prop:limit-empty-iff-some-empty}: towards this, we recall some notions from model theory, prove the intermediate Proposition \ref{prop:compact-in-extension}, and describe how to obtain an infinite system of algebraic equations from a differential ideal. In Section \ref{section:theorem} we state and give a proof for the first equality of the theorem, i.e.\ that the set of solutions to the tropicalization of a differential equation and the set of tropicalizations of the classical solutions coincide. In Section \ref{section:grobner} we extend Gröbner theory to the differential nontrivially valued case and conclude adding the second equality to the fundamental theorem, stating that the two aformentioned sets coincide with the set of weight vector whose initial ideals do not contain monomials. 
Finally, in Section \ref{section:raggi} we introduce the radius of convergence of a tropical power series and discuss Corollary \ref{corolla:intro}, giving an explicit example of computation of tropical solutions to a differential equation. We conclude with a table of notations to help the reader navigate the paper.

\subsection*{Acknowledgements}
We would like to thank Jeff Giansiracusa for many helpful discussions and comments and for carefully reading an earlier draft of this paper, Marie-Charlotte Brandenburg for providing a proof of Lemma \ref{lem:finite-vertices}, and Ángel David Ríos Ortiz, Alejandro Vargas, Fuensanta Aroca and Benjamin Smith for helpful discussions and comments during the preparation of this paper.
S.M. is thankful to MPI MiS Leipzig for the excellent working conditions during the preparation of this paper and for its stimulating mathematical environment.
S.M. was partially supported by the Wallenberg AI, Autonomous Systems and Software Program (WASP) funded by the Knut and Alice Wallenberg Foundation.
F.G. is grateful to Amador Martin-Pizarro for bringing this problem to his attention.
F.G. was partially supported by the program GeoMod ANR-19-CE40-0022-01 (ANR-DFG) and by the Italian research project PRIN 2022 ``Models, sets and classifications'', Prot. 2022TECZJA.

	\section{Tropical differential equations} \label{section:tropdiffeq}

In order for this work to be as self-contained as possible, in this section we give a presentation of the framework in which tropical differential equations (also shortened as TDEs in the following) are considered and how tropicalization is defined for algebraic differential equations with coefficients in a differential ring $R$ and for their solutions. 

\subsection*{Algebraic differential equations}
From here on, $m$ will be a nonzero natural number, and $[m]$ will denote the set $\{1, \dots , m\}$. A \emph{partial differential ring} with $m$ commuting differentials $(R, \{ d_{R,i}\}_{i\in [m]})$ is a ring equipped with $m$ pairwise commuting additive maps 
\[
	d_{R,1}, \dots , d_{R,m} \co R \rightarrow R
\]
that satisfy the Leibniz rule. Given an element $x$ in a partial differential ring $(R, \{d_{R,i}\}_{i \in [m]})$ and a vector $\mathbf j \in \N^m$ we will denote by $d^{\mathbf{j}}x$ the element 
\[
d_{R,1}^{j_1}(d_{R,2}^{j_2}( \dots (d_{R,m}^{j_m}(x)) \dots )).
\]
In the present work we will only deal with differential rings of the form 
\[
	(\pkk m, \{d/dt_i\}_{i \in [m]})
\]
for $K$ a field. In the following we denote $\pkk m$ by $R^K_m$ (and by $R_m$ if the field of coefficients is clear).

A \emph{differential algebra} $(A, \{d_{A,i}\}_{i \in [m]})$ over $R_m$ is a partial differential ring with $m$ commuting differentials, equipped with a homomorphism $f \co R_m \rightarrow A$ such that 
\[
    f\left (\frac{d}{dt_i}(r) \right) = d_{A,i}(f(r))
\]
for all $i \in [m]$ and for all $r \in R_m$.
The main example of differential algebra over $R_m$ is the differential ring of \emph{differential polynomials} in $n$ variables over $R_m$, as introduced by Ritt in \cite{ritt}. It is defined as follows:
\begin{equation}\label{eq:ritt-polynomials}
\diff {R_m} n := (R_m[x_i^{(\mathbf j)} \mid i \in [n] ; \mathbf j \in \N^m], \{d_i\}_{i \in [m]})
\end{equation}
where the differential $d_i$ extends $d/dt_i$ for each $i= 1, \dots , m$ and acts on the variables by sending $x_k^{(\mathbf j)}$ to $x_k^{(\mathbf j + e_i)}$ for all $k \in [n]$, where $e_i$ is the vector in $\N^m$ whose only nonzero entry is the $i$-th, which is 1. We will often write $x_i$ for $x_i^{(\mathbf 0)}$. When considering $R_m[x_i^{(\mathbf j)} \mid i \in [n] ; \mathbf j \in \N^m]$ just as an $R_m$-algebra forgetting about its differential structure, we will denote it as  $\basic {R_m} n$: the reason for introducing this notation will be clear later, when we will introduce Ritt algebras over partial differential semirings.

An ideal $I$ of $R_m$ is \emph{differential} if it is closed under differentiation.
Given a differential algebra $A$ over $R_m$, a solution in $A$ to the differential equation $f \in \diff {R_m} n$ is a $n$-tuple $a=(a_1, \dots, a_n) \in A^n$ such that plugging in $d_A^{\mathbf j}(a_i)$ for $x_i^{(\mathbf j)}$ in $f$ results in $0$. We will denote as $\textup{Sol}_A(E)$ the set of solutions to a system of differential equations $E$ in $A$.

\subsection*{Semirings and valuations}
Let us now introduce the main objects in the theory of TDEs, in order to state their rigorous definition and make clear how to tropicalize an algebraic differential equation to get a tropical one.
Among the several algebraic foundations for tropical geometry, such as  hyperfields
\cite{Viro, Baker-Bowler,Lorscheid-hyperfield,Jun-hyperfields, james, jellscheidereryu}, Lorscheid's blueprints
\cite{Lorscheid-blueprints1,Lorscheid-scheme-theoretic}, and idempotent semirings
\cite{eqtrop,univtrop,GG3,Maclagan-Rincon-1,Maclagan-Rincon-2,Noah-module-theoretic, Yaghmayi}, following the presentation developed in \cite{gianmereta} and \cite{tesimereta}, we will use the language of (differential) idempotent semirings. 

Firstly, let us recap the definition of semiring and of valuation and give some relevant examples: a \emph{semiring} $(S, \oplus, \odot)$ is an algebraic structure satisfying the same axioms as a ring, except for the existence of additive inverses. A semiring is \emph{idempotent} if $a \oplus a = a$ for all elements $a \in S$. When $S$ is idempotent the operation $\oplus $ induces a
canonical partial order on $S$ defined by $a \preceq b$ if $a \oplus b = b$. We will say that a sum in an idempotent semiring \emph{tropically vanishes} if the result of the sum after removing any one of the summands does not change.
\begin{ex}
\begin{enumerate}
\item The set $\B = \{0,\infty\}$ equipped with operations $\oplus := \min$ and $\odot := +$ is an idempotent semiring, called the idempotent semiring of Boolean numbers;
\item The set $\T = \R \cup \{\infty\}$ equipped with operations $\oplus := \min$ and $\odot: = +$ is an idempotent semiring, called the idempotent semiring of tropical numbers. It contains $\B$ as a subsemiring;
\item For every $n \ge 1$ let $\T^{(n)} = \R^n \cup \{\infty\}$ with $\oplus$ given by lexicographic minimum and $\odot := +$. It is an idempotent semiring.
\end{enumerate}
\end{ex}
\begin{oss}
    Note that the partial order $\preceq$ induced by $\oplus := \min$ in the examples above is the reverse of the order $\le$ with respect to which we take the minimum.
\end{oss}

In order to introduce the main examples of semirings we will work with, we need to recall some properties of cones in $\R^m$:
\begin{defin}
A \emph{rational polyhedral cone} $\tau \subset \R^m$ is a subset of the form:
\[
    \{\lambda_1 v_1 + \dots + \lambda_n v_n \mid \lambda_i \in \R_{\ge 0}\text{ for all } i \in [n]\}
\]
for vectors $v_1, \dots , v_n \in \Z^m$.  A cone $\tau$ is \emph{strongly convex} if the only vector subspace contained in $\tau$ is the trivial subspace of $\R^m$. Finally, we say that a cone has dimension $r$ if the smallest subspace of $\R^m$ containing $\tau$ has dimension $r$.
From now on, the word \emph{cone} is intended to mean \emph{strongly convex rational polyhedral cone, of dimension $m$}.  
\end{defin}
Let us introduce some notation for the following: 
\begin{itemize}
\item given a subsemiring $S$ of $\T$, the support of a function $f \co \Z^m \rightarrow S$ is 
\[
    \supp (f) := \{\jj \in \Z^m \mid f(\jj) \neq \infty\} \subseteq \Z^m.
\]
Analogously, the support of a formal series $A=\sum_{\jj \in \Z^m}a_\jj t^\jj$ with coefficients in $S$ is 
\[
    \supp (A) := \{\jj \in \Z^m \mid a_\jj \neq \infty \} \subseteq \Z^m;
\]
\item given $Z_1, Z_2 \subseteq \R^m$, we denote by $Z_1 + Z_2$ their Minkowski sum, i.e.\ the set 
\[
    \{z_1+z_2 \mid z_1 \in Z_1, z_2 \in Z_2\} \subseteq \R^m.
\] 
\end{itemize}
\begin{defin}
    Given a cone $\tau \subset \R^m$, and a subset $Z$ of $\R^m$, the \emph{convex hull} $C_\tau(Z)$ of $Z$ with respect to $\tau$ is the smallest convex set in $\R^m$ containing $Z + \tau$.
    An element $z \in Z$ is a \emph{vertex} of $C_\tau(Z)$ if $C_\tau(Z \setminus \{z\}) \subsetneq C_\tau(Z)$. Let us denote the set of vertices of $C_\tau(Z)$ as $\textup{Vert}_\tau(Z)$.
\end{defin}
We now introduce the main examples of semirings in our work:
\begin{defin}\label{def:conv-boolean-semirings}
    Given a cone $\tau \subset \R^m$ and a subsemiring $S$ of $\T$, let $\textup{Conv}(\tau,S)$ be 
    \[
        \textup{Conv}(\tau,S) := 
        \{(C_\tau(Z), \alpha_Z \co \Z^m \cap \tau  \rightarrow S) \mid Z \subset \Z^m \cap \tau \text{ finite, supp}(\alpha_Z) =\textup{Vert}_\tau(Z)\}.
    \] 
    
    The set $\textup{Conv}(\tau,S)$ is an idempotent semiring with the following operations: 
    \begin{itemize}
    \item for $Z_1, Z_2$ finite subsets of $\Z^m \cap \tau$, let 
    \[
        (C_\tau(Z_1),\alpha_{Z_1})\oplus(C_\tau(Z_2), \alpha_{Z_2}) := (C_\tau(Z_1 \cup Z_2), \alpha_{Z_1 \cup Z_2})
    \]
    with
    \[
        \alpha_{Z_1 \cup Z_2}(\jj) := \begin{cases}
        \alpha_{Z_1}(\jj) \oplus \alpha_{Z_2}(\jj) & \text{if } \jj \in \textup{Vert}_\tau(Z_1 \cup Z_2)  
        \\
        \infty & \text{if } \jj \notin \textup{Vert}_\tau(Z_1 \cup Z_2)  
        \end{cases};
    \]
    \item for $Z_1, Z_2$ finite subsets of $\Z^m \cap \tau$, let 
    \[
        (C_\tau(Z_1),\alpha_{Z_1})\odot(C_\tau(Z_2), \alpha_{Z_2}) = (C_\tau(Z_1 + Z_2), \alpha_{Z_1 + Z_2})
    \]
    with
    \[
    \alpha_{Z_1 + Z_2}(\jj) := \begin{cases}
    \alpha_{Z_1}(\jj) \odot \alpha_{Z_2}(\jj) & \text{if } \jj \in \textup{Vert}_\tau(Z_1 + Z_2)  
    \\
    \infty & \text{if } \jj \notin \textup{Vert}_\tau(Z_1 + Z_2)  
        \end{cases}.
    \]
    \end{itemize}
    With these operations, we have $0_{\textup{Conv}(\tau, S)} = (\varnothing, \infty)$ and $1_{\textup{Conv}(\tau, S)}=(\tau=C_\tau(\{0_{\R^m}\}),\alpha_{\{0_{\R^m}\}})$, with 
    \[
    \alpha_{\{0_{\R^m}\}}(\jj) := \begin{cases}
    0 & \text{if } \jj = 0_{\R^m}
    \\
    \infty & \text{otherwise} . 
        \end{cases}
    \]
    The semiring $\textup{Conv}(\tau,S)$ is an $S$-algebra via the injective morphism mapping each element $s \in S$ to $(\tau,s \odot_S \alpha_{\{0_{\R^m}\}})$.
\end{defin}
\begin{oss}\label{oss:representation-of-conv}
Let $\tau \subset \R^m$ be a cone and let $S$ as above. Given $C := (C_{\tau}(Z), \alpha_Z) \in \textup{Conv}(\tau,S)$, let $\ell$ be the number of vertices of $C_{\tau}(Z)$. We can represent $C$ as 
\begin{equation}\label{eq:repr-convex-hull}
    C = ((\jj_{1}; \alpha_1) ,\dots ,(\jj_\ell;\alpha_\ell) )
\end{equation}
where $\textup{Vert}_\tau(Z) = \{\jj_{1}, \dots ,\jj_\ell\}$ is the set of vertices of $C$ and $\alpha_i := \alpha_Z(\jj_i) \in S \setminus\{0_S\}$. We will refer to $\alpha_i$ as the \emph{weight} of $\jj_i$. The representation \ref{eq:repr-convex-hull} is unique up to permutation. When $C$ has only one vertex and one weight, we will simply write $C=(\jj; \alpha)$. 
\end{oss}
When $\tau = \posR$ in the construction above, we will denote $\textup{Conv}(\tau,S)$ as $\textup{Conv}_m(S)$. The idempotent semiring $\textup{Conv}_m(\B)$ appears already in \cite{cotterill} denoted as $V\B(\mathbf t)$, and, with different notation, in \cite{sebastian}.
\begin{oss}\label{oss:Conv-m=1}
Notice that for $m=1$ we have the following isomorphisms of idempotent semirings: 
\[
    \textup{Conv}_1(\B) \cong \N \cup \{\infty\} \subset \T  \quad \quad \quad \textup{Conv}_1(\T) \cong (\N \times \R) \cup \{\infty\} \subset \T^{(2)}.
\]
In the first case this is implemented by sending the convex hull of a finite set $Z \subset \tau \cap \Z=\N$ to its only vertex (i.e.\ the minimum of $Z$) and $\varnothing$ to $\infty$. In the second case it is by sending a convex hull with vertex $z$ and vertex weight $\alpha$ to the element $(z, \alpha) \in \T^{(2)}$ and  $\varnothing$ to $\infty$. 

Furthermore, for every cone $\tau$ and idempotent semiring $S$, there is an homomorphism of semirings
    \begin{equation}\label{eq:sigma-0}
        \sigma_0 \co \textup{Conv}(\tau,S) \rightarrow \textup{Conv}(\tau,\B)
    \end{equation}
induced by the canonical morphism $S \rightarrow \B$, i.e.\ by forgetting the weights. For an element $C \in \textup{Conv}(\tau,S)$ we will refer to $\sigma_0(C)$ as the \emph{support} of $C$. 
\end{oss}
\begin{ex}
We give an example of sum and product of two elements of $\textup{Conv}_m(\T)$. Let $C_1, C_2 \in \textup{Conv}_2(\T)$ be defined as
\[
C_1 = (((1,2);0) ,((5,0);2) )
\]
\[
C_2 = ( ((0,4);1), ((1,2);1), ((3,1);1) ).
\]
We can represent them pictorially as follows:
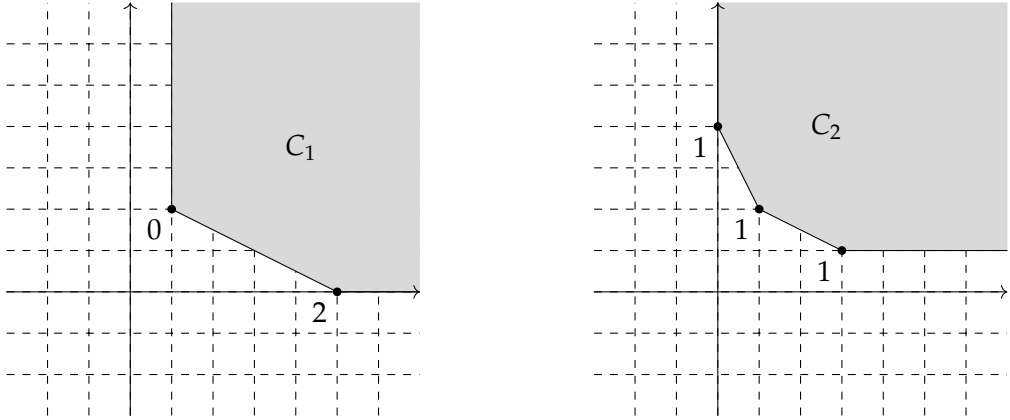
\begin{center}
\hspace{-2.5cm}
\begin{minipage}{.2\textwidth}
\begin{tikzpicture}[scale=0.45]
\draw (-2.99,-2.99)[dashed, opacity=0.5] grid (6.99,6.99);
\fill[fill=gray!30, opacity=0.75] (1,7)--(1,2)--(5,0)--(7,0)--(7,7);
\draw[->] (-3,0) -- (7,0);
\draw [->] (0,-3) -- (0,7);
\draw (1,7) -- (1,2);
\draw (1,2)--(5,0);
\draw (5,0)--(7,0);
\draw (3.5,3.5) node[right] {$C_1$};
\draw (5,0) node[below left] {2};
\draw (1,2) node[below left] {0};
\fill (5,0) circle[radius=3pt];
\fill (1,2) circle[radius=3pt];
\end{tikzpicture}
\end{minipage} \hspace{4.5cm}
\begin{minipage}{.2\textwidth}
\begin{tikzpicture}[scale=0.45]
\draw (-2.99,-2.99)[dashed, opacity=0.5] grid (6.99,6.99);
\fill[fill=gray!30, opacity=0.75] (0,7)--(0,4)--(1,2)--(1,2)--(3,1)--(7,1)--(7,7);
\draw[->] (-3,0) -- (7,0);
\draw [->] (0,-3) -- (0,7);
\draw (0,7) -- (0,4);
\draw (0,4) -- (1,2);
\draw (3,1)--(1,2);
\draw (3,1)--(7,1);
\draw (2,4) node[right] {$C_2$};
\draw (0,4) node[below left] {1};
\draw (1,2) node[below left] {1};
\draw (3,1) node[below left] {1};
\fill (0,4) circle[radius=3pt];
\fill (1,2) circle[radius=3pt];
\fill (3,1) circle[radius=3pt];
\end{tikzpicture}
\end{minipage}
\end{center}
Their sum and product are respectively:
\[
C_1 \oplus C_2 = (((0,4);1) ,((1,2);0), ((5,0);0) ))
\]
\[
C_1 \odot C_2 = ( ((1,6);1),((2,4);1), ((8,1);3) ))
\]
and are represented below:
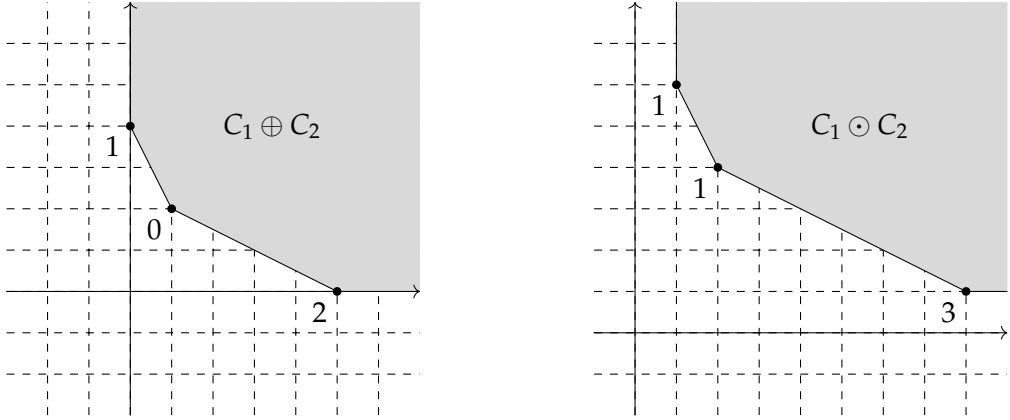
\begin{center}
\hspace{-2.5cm}
\begin{minipage}{.2\textwidth}
\begin{tikzpicture}[scale=0.45]
\draw (-2.99,-2.99)[dashed, opacity=0.5] grid (6.99,6.99);
\fill[fill=gray!30, opacity=0.75] (0,7)--(0,4)--(1,2)--(5,0)--(7,0)--(7,7);
\draw[->] (-3,0) -- (7,0);
\draw [->] (0,-3) -- (0,7);
\draw (0,7) -- (0,4);
\draw (0,4)--(1,2);
\draw (1,2)--(5,0);
\draw (2,4) node[right] {$C_1 \oplus C_2$};
\draw (0,4) node[below left] {1};
\draw (1,2) node[below left] {0};
\draw (5,0) node[below left] {2};
\fill (0,4) circle[radius=3pt];
\fill (1,2) circle[radius=3pt];
\fill (5,0) circle[radius=3pt];
\end{tikzpicture}
\end{minipage} \hspace{4.5cm}
\begin{minipage}{.2\textwidth}
\begin{tikzpicture}[scale=0.45]
\draw (-2.99,-2.99)[dashed, opacity=0.5] grid (6.99,6.99);
\fill[fill=gray!30, opacity=0.75] (-1,7)--(-1,5)--(0,3)--(4,1)--(6,0)--(7,0)--(7,7);
\draw[->] (-3,-1) -- (7,-1);
\draw [->] (-2,-3) -- (-2,7);
\draw (-1,7) -- (-1,5);
\draw (0,3) -- (-1,5);
\draw (0,3)--(4,1);
\draw (4,1)--(6,0);
\draw (7,0)--(6,0);
\draw (2,4) node[right] {$C_1\odot C_2$};
\draw (-1,5) node[below left] {1};
\draw (0,3) node[below left] {1};
\draw (6,0) node[below left] {3};
\fill (-1,5) circle[radius=3pt];
\fill (0,3) circle[radius=3pt];
\fill (6,0) circle[radius=3pt];
\end{tikzpicture}
\end{minipage}
\end{center}
\end{ex}
We adopt the following generalised definition of valuation, as introduced in \cite{eqtrop, notejeff} and used in \cite{gianmereta}: a \emph{valuation} on a ring $R$ is given by an idempotent semiring $S$ and a map $v\co R \to S$ satisfying
\begin{enumerate}
\item	$v(0)  = 0_S$,
\item	$v(1) = v(-1) = 1_S$,
\item	$v(ab) = v(a) \odot v(b)$,
\item	$v(a+b) \oplus v(a) \oplus v(b)$ tropically vanishes (i.e.\ $ v(a+b) \preceq v(a) \oplus v(b)$).
\end{enumerate}
Often we will refer to the datum $(R,v)$ of a ring with a valuation as a \emph{valued ring}. 
\begin{oss}\label{oss:rank2valuation}
    Notice that a valuation $v \co R \rightarrow \T$ can be extended to a valuation $R\{\!\{t\}\!\} \rightarrow \T^{(2)}$ on the ring $R\{\!\{t\}\!\}$ of Puiseux series (or the subrings of formal Laurent series, power series or polynomials) by the formula
\[
	a_0 t^{n_0} + \cdots \mapsto (n_0, v(a_0)) \in \T^{(2)}.
\]
\end{oss}
In order to introduce the main examples of valuations we will be working with, we need the following definition and lemma. 
\begin{defin}
    Let $\tau$ be a cone as above and let $S$ be an idempotent semiring. Define $S[\![\tau]\!]$ as follows:
    \[
    S[\![\tau]\!] := \left \{ A = \sum_{\jj \in  \Z^m} s_\jj t^\jj \mid \supp(A) \subseteq  \tau  \right \}.
    \]
    This set, equipped with usual sum and product, is the idempotent semiring of multivariate power series supported on $\tau$. When $S$ is a field, this was introduced (with fractional coefficients) for example in \cite{arocacanojung, arocarond}. 
    When $\tau = (\R_{\ge 0})^m$, we will denote $S[\![\tau]\!]$ as $S[\![t_1, \dots t_m]\!]$.
\end{defin}
Notice that restricting the support of a power series to a cone is necessary to make the usual product of power series a well-defined operation.
\begin{lem}\label{lem:finite-vertices}
    Let $\tau \subset \R^m$ be a cone and $Z \subset \tau \cap \Z^m$. Then $\textup{Vert}_\tau(Z)$ is finite.
\end{lem}
\begin{proof}
    For an element $z \in Z$ we say that $z$ is a \emph{semivertex} if $(Z \setminus \{z\})+ \tau \subsetneq Z + \tau$. Clearly every element in $\textup{Vert}_\tau(Z)$ is a semivertex. Find nonzero (without loss of generality, rational) vectors $v_1, \dots , v_m \in \R^m$ such that $\tau$ can be written as $\{x \in \R^m \mid \langle v_i, x \rangle \ge 0 \text{ for each } i \in [m]\}$. Then for all $z \in \R^m$,
    \[
        z + \tau =\{x \mid x-z \in \tau\} = \{x \mid \langle v_i, x \rangle \ge \langle v_i, z \rangle \text{ for each } i \in [m]\}.
    \]
    For $i \in [m]$, let $\alpha_i = \min_{z \in Z}\langle v_i, z\rangle \ge 0$. This minimum exists since $\tau$ is rational, so the normal vectors $v_i$ are rational and after scaling them, we can assume that they are integral. Every lattice point $x \in \tau$ has $\langle x,v_i \rangle \ge 0$, and since both $x$ and $v_i$ are integral, so is their inner product. As a decreasing sequence of non negative integers stabilizes, the minimum exists. 

    Let $D :=\{x \in \R^m \mid \langle v_i,x\rangle \ge \alpha_i \text{ for each } i \in [m]\}$. It is a polyhedron with recession cone $\tau$. By construction we have $\tau \supseteq D \supseteq Z + \tau$, and the convexity of $D$ implies $D \supseteq C_\tau(Z)$. For $i \in [m]$, let $z_i$ be such that $\langle v_i, z_i \rangle = \alpha_i$ and $\beta_i := \max_{j \in [m]}\langle v_i, z_j\rangle$.
    Let $E := \{x \in \R^m \mid \alpha_i \le \langle  v_i, x\rangle \le \beta_i \text { for each } i \in [m]\}$. This is a bounded subset of $D$. 
    
    We prove the statement by showing that every semivertex is contained in $E$. Note that if $x \in \bigcup_{j \in [m]} (z_j + \tau)$ and $x \neq z_j$ for all $j$, then $x$ is not a semivertex. Hence, semivertices are contained in the set
    \[
        \left ( D \setminus \bigcup_{j=1}^m (z_j + \tau) \right ) \cup \{z_1, \dots , z_m\}.
    \]
    We now show that the set above is contained in $E$:
    \begin{align*}
        D \supseteq E \cup \bigcup_{j=1} ^m (z_j + \tau) & = E \cup \bigcup_{j=1}^m\{x \in \R^m \mid \langle v_i , x \rangle \ge \langle v_i, z_j \rangle \text{ for all } i \in [m]\} = \\
        & =
        E \cup \{x \in \R^m \mid \exists j \in [m] \text{ s.t. } \langle v_i , x \rangle \ge \langle v_i, z_j \rangle \text{ for all } i \in [m]\} \supseteq 
        \\
        &  \supseteq 
        E \cup \{x \in \R^m \mid \langle v_i , x \rangle \ge \beta_i \text{ for all } i \in [m]\} =
        \\
        & = D.
    \end{align*}
    So we get
    \[
        E \supseteq \overline{ \left (D \setminus \bigcup_{j=1}^m (z_j + \tau) \right ) } \supseteq \left ( D \setminus \bigcup_{j=1}^m (z_i + \tau) \right ) \cup \{z_1, \dots , z_m\}.
    \]
    Thus all the semivertices lie in $E$. This implies that $\textup{Vert}_\tau(Z) \subseteq E$ and since $E$ is bounded, $\textup{Vert}_\tau(Z)$ is finite.
\end{proof}
\begin{oss}
    The Lemma above holds with analogous proof also allowing fractional exponents, since the definition of fractional Laurent series, as it can be found in \cite{arocacanojung}, implies the existence of a least common multiple of the denominators of the exponents appearing in each series.
\end{oss}
\begin{oss} \label{oss:conv-as-quotient}
 Given a cone $\tau$ and an idempotent semiring $S$, it follows from Lemma \ref{lem:finite-vertices} that the idempotent semiring $\textup{Conv}(\tau,S)$ of Definition \ref{def:conv-boolean-semirings} can be realised as the quotient of $S[\![\tau]\!]$ by the semiring congruence
\begin{center}
\begin{tabular}{ c c c }
 & & $C_\tau(\supp(A)) = C_\tau(\supp(B)) $ \\ 
 $A = \sum a_\jj t^\jj  \sim  B =\sum b_\jj t^\jj $ & $\iff$ & and \\  
 &  & $a_\jj = b_\jj \text{ for all } \jj \in \textup{Vert}_\tau(\supp(A))$   
\end{tabular}
\end{center}
When $S=K$ is a field, given $A \in \pkcone \tau$, we have that the convex hull of the support of $A$ with respect to $\tau$ has finitely many vertices. This allow us to define the following valuations:
\end{oss}
\begin{ex}\label{ex:conv-vals}
\begin{enumerate}
    \item Given a cone $\tau \subset \R^m$ and a field $K$, the map 
            \[
	           w_\tau \co  \pkcone \tau  \to \textup{Conv}(\tau,\B)
		\]
    sending an element $A \in \pkcone \tau $ to $C_\tau(\supp(A))$  is a valuation. As noted in Remark \ref{oss:Conv-m=1}, when $m=1$, $\textup{Conv}_1(\B) \subset \T$ and it is easy to check that the valuation $w$ becomes the $t$-adic valuation $R_1\rightarrow \T$ sending a univariate power series to its leading order.
    \item  Given a cone $\tau \subset \R^m$, let $v_K \co K \rightarrow \T$ be a valued field. With a similar procedure to that illustrated in Remark \ref{oss:rank2valuation}, we build a refined version of the valuation $w_\tau$ above. The map
	\[ 
		v_\tau \co \pkcone \tau \to \textup{Conv}(\tau,\T)
	\] 
    given by sending $A \in \pkcone \tau$ to $C_\tau(\supp(A))$, with vertices weighted by the valuation $v_K(a_{\mathbf{n}})$ of the coefficients of the corresponding terms, is a valuation.
    Analogously as for the valuation $w$ of point (1), for $m=1$, $\textup{Conv}_1(\T) \subset \T^{(2)}$ and we retrieve the valuation $R_1\rightarrow \T^{(2)}$ obtained by refining the $t$-adic valuation as in Remark \ref{oss:rank2valuation}. 
\end{enumerate}

\end{ex}

\begin{oss}
    For every cone $\tau$, the valuation $w_\tau$ above is surjective. On the other hand, $v_\tau$ is surjective if and only if $v_K$ is a surjective valuation (or if we restrict the weights to the value group of $v_K$). Furthermore, notice that for a valued field $v_K \co K \rightarrow \T$ the diagram
    \[
        \begin{tikzcd}
            & \pkcone \tau  \arrow[dr, "w_\tau"] \arrow[dl, swap, "v_\tau"]&  \\
             \textup{Conv}(\tau,\T) \arrow[rr, "\sigma_0"]  & & \textup{Conv}(\tau, \B) 
        \end{tikzcd}
    \]
    commutes for every cone $\tau$.
\end{oss}

\begin{defin}
    Given a totally ordered commutative monoid $(G, +, \preceq)$, let $\widetilde G$ be the idempotent semiring $(G \cup \{\infty\}, \min_{\preceq}, +)$.
\end{defin}

We introduce here another semiring and a valuation that will turn out to be useful in Section \ref{section:grobner}. Let $\omega \in \R^m$ with rationally independent coordinates. Then the relation
    \[
        v_1 \preceq_\omega v_2 \iff \langle v_1 , \omega \rangle \le \langle v_2, \omega \rangle,
    \]
is a total order on any additive submonoid $G$ of $\Q^m$. 
\begin{defin} \label{def:lead-semiring}
Given a subsemiring $S$ of $\T$ and a submonoid $G$ of $\Q^m$, endow the commutative monoid $G \times (S \setminus\{0_S\})$ with the total order $\preceq$ defined as the lexicographic order induced by $\preceq_\omega$ and $\preceq_S$. Set 
\[
\textup{Lead}_\omega(G,S):=\widetilde{G \times (S \setminus\{0_S\})}.
\]
Notice that $\textup{Lead}_\omega(G,S)$ is an $S$-algebra.
\end{defin}

\begin{ex} \label{ex:omega-val-tau}
Given $\omega \in \R^m$ with rationally independent coordinates, let $H_\omega$ be the halfspace $\{x \in \R^m \mid \langle x, \omega \rangle \ge 0 \}$. Given a cone $\tau \subset H_\omega$, the vector $\omega$ induces the total  order $\preceq_\omega$ on $\tau \cap \Z^m$ as above.
Let $v_K \co K \rightarrow \T$ be a valued field, the map 
    \[
        v_{\tau,\omega} \co \pkcone \tau \rightarrow \textup{Lead}_\omega(\tau \cap \Z^m,\T)
    \]
    defined as:
    \[
        A=\sum_{\jj \in \tau \cap \Z^m} a_\jj t^\jj \mapsto (\overline \jj, v_K (a_{\overline \jj})),
    \]
    where $\overline \jj= \min_{\preceq_\omega} \supp(A)$, is a valuation. When $m=1$, this is the valuation $R_1 \rightarrow \T^{(2)}$ refining the $t$-adic valuation as in Remark \ref{oss:rank2valuation}.
\end{ex}
\subsection*{Differential semirings}
Given an idempotent semiring $S$, an
additive map $d\co S \to S$ is said to be a \emph{tropical differential} if it satisfies the tropical Leibniz relations:  for any two elements  $x,y \in S$ the expression
\[
d(xy) \oplus xd(y) \oplus yd(x)
\]
tropically vanishes. Note that we can view the tropical Leibniz relations as the tropicalization of the usual Leibniz rule. Furthermore, as $S$ is an idempotent semiring, if an additive map $d\co S \to S$ satisfies the usual Leibniz rule, it satisfies the tropical Leibniz relations. In this case we will say that $d$ is a \emph{strict} tropical differential.

\begin{defin}
A \emph{partial differential idempotent semiring} $(S, \{ d_{S,i}\}_{i \in [m]} \})$ with $m$ commuting differentials is an idempotent semiring $S$ equipped with $m$ pairwise commuting tropical differentials $d_{S,i} \co S \rightarrow S$. Analogously as for rings, a homomorphism of differential semiring is a homomorphism of semirings commuting with the differentials.

We will often drop the adjectives \emph{partial} and \emph{idempotent} in the following.
\end{defin}

The two main examples of differential semiring we will work with are the following:
\begin{ex}\label{ex:diffsemirings}
\begin{enumerate}
\item Consider the idempotent semiring $\pbb m$ of multivariate formal power series with boolean coefficients. Endowing it with the differentials $\{d/dt_i\}_{i \in [m]}$ defined by
\[
\frac{d}{dt_i}(t_j^n) = \begin{cases} 
t_i^{n-1} & n \geq 1, i=j\\
\infty & \text{otherwise}
\end{cases}
\]
and extended by Leibniz rule we obtain a strict partial differential semiring. In the following, we denote $(\pbb m, \{d/dt_i\}_{i \in [m]})$ as $\B_m$.
\item Consider the idempotent semiring of multivariate tropical power series $\ptt m$. Given a nontrivial valuation $v\co \mathbb{N} \to \T$, let $(d/dt_i)_v$ be defined as 
\[
\left (\frac{d}{dt_i} \right)_v(t_j^n) = \begin{cases} 
v(n)t_i^{n-1} & n \geq 1, i=j\\
\infty & \text{otherwise},
\end{cases}
\]
then $(d/dt_i)_v$ is a tropical differential for every $i$. Denote as $\T_m$ the semiring $\ptt m$ and as $\T_{m,v}$ the (non-strict) differential semiring $(\ptt m, \{(d/dt_i)_v\}_{i \in [m]})$. 
\end{enumerate}
\end{ex}

The definition of the differential semiring 
 \[
    (\diff S n , \{d_i\}_{i \in [m]})
 \]
of differential polynomials over $(S,\{d_{S,i}\}_{i \in [m]})$ is quite convoluted and we are not going to recall it here. All the details can be found in \cite[Section 3.4]{gianmereta} or \cite[Chapter 8]{tesimereta}.
 It satisfies the desired universal property and it gives back Ritt's construction as in \ref{eq:ritt-polynomials} when  $(S,\{d_{S,i}\}_{i \in [m]})$ is strict (in particular for differential rings and $\B_m$). In the general case the $S$-algebra inclusion
\[
    \basic S n \subseteq \diff S n .
\]
holds, and can be strict.

\subsection*{Tropical pairs}
The central objects in the theory of TDEs are the so-called tropical pairs, playing in this theory an analogous role to that of rings and algebras in the theory of polynomial equations: 
\begin{defin}
A \emph{(partial) tropical pair} $\mathbf{S}$ with $m$ differentials is a partial differential semiring $(S_1,\{d_{S_1,i}\}_{i \in [m]})$ with $m$ commuting differentials and a homomorphism of idempotent semirings $\pi\co S_1 \to S_0$. We will say that $\mathbf S$ is strict if $S_1$ is a strict differential semiring. A morphism of pairs $\sigma$ from $(S_1\to S_0)$ to $(T_1
\to T_0)$ is a commutative diagram of idempotent semirings
\[
\begin{tikzcd}
S_1 \arrow[d] \arrow[r, "\sigma_1"] & T_1 \arrow[d] \\
S_0 \arrow[r, swap, "\sigma_0"]           & T_0          
\end{tikzcd}
\]
in which the upper horizontal arrow $\sigma_1$ is a morphism of partial differential idempotent semirings.
\end{defin}
These maps can be thought of as sending elements of $S_1$ to their leading orders. This will actually be the case for the instances considered in the remainder of this work. We introduce them in the following example:

\begin{ex}\label{ex:pairs}
\begin{enumerate}
    \item Consider $(\B_1, d)$ as in point (1) of Example \ref{ex:diffsemirings} for $m=1$. The homomorphism $\Psi \co \B_1 \to \T$ defined by $t^n \mapsto n$ is a tropical pair.
    \item Consider $ \Phi\co \pt_v\to \T^{(2)}$, where the source has any of the differentials from point (2) of Example \ref{ex:diffsemirings} and the morphism $\Phi$ is given by
    \[
        (a_{n_0} t^{n_0} + a_{n_1} t^{n_1} + \cdots) \mapsto (n_0, a_{n_0}).
    \] 
    It is a tropical pair.
    \item For every $m$, let $\B_m$ as in Example \ref{ex:diffsemirings} and consider the partial tropical pair
            \[
	           \Psi_m \co  \B_m \to \textup{Conv}_m(\B).
		\]
    sending an element $B \in \pbb r$ to the convex hull of its support. For $m=1$, analogously to point (1) of Example \ref{ex:conv-vals}, $\Psi_1$ is the pair $\Psi$ of point (1). We will denote this pair by $\mathbf T _m$.
    \item Consider the partial tropical pair
		\[ 
		\Phi_m \co \T_{m,v} \to \textup{Conv}_m(\T),
		\] 
    where the source has differentials as in Example \ref{ex:diffsemirings} for a valuation $v$. The morphism $\Phi_m$ is given by sending $B \in \T_{m,v}$ to the convex hull of its support with vertices weighted by the coefficients of the corresponding leading terms.
    When $m=1$, $\Phi_1$ is the pair $\Phi$ of point (2). We will denote this pair as $\mathbf S _m$.
    \end{enumerate}
\end{ex}

All pairs appearing in Example \ref{ex:pairs} are \emph{reduced} (see \cite[Section 4]{gianmereta} for definition and properties of reduced pairs and reduction). Since all the pairs we work with are reduced, from now on every instance of the word "pair" has to be intended as "\emph{reduced} pair". 

Given a pair $\mathbf{S}$, the category of $\mathbf{S}$-algebras is the category of pairs $\mathbf{T}$ under $\mathbf{S}$. A typical example of an $\mathbf S$-algebra is $\diff {\mathbf S} n$, the algebra of differential polynomials in $n$ variables over $\mathbf S$. 
The pair $\diff {\mathbf{S}} n$ is obtained as the reduction of the pair
    \[
        \diff {S_1}n \rightarrow (S_0 | S_1)\{x_1, \ldots, x_n\}.
    \]
where $(S_0 | S_1)\{x_1, \ldots, x_n\}$ is obtained as the pushout of a certain diagram: we refer again to \cite[Section 4]{gianmereta} for the details of this construction. For our scope, it is enough to say that the inclusion of $S_0$-algebras
\[
    \basic {S_0} {n} \subseteq (S_0 | S_1)\{x_1, \ldots, x_n\}.
\]
holds.
The $\mathbf{S}$-algebra $\diff {\mathbf S} n$ satisfies a universal property similar to that of Ritt's polynomials (see \cite[Proposition 4.1.4]{gianmereta}).

\subsection*{TDEs, differential enhancements and tropicalization}
We can finally state formally the definition of tropical differential equation that we will use, and how solutions to such an equation look like:
\begin{defin}\label{defin:tropical-solutions}
\begin{enumerate}
    \item A \emph{tropical differential equation} (TDE) $f$ over the pair $\mathbf S = (S_1 \stackrel{\Phi}{\to} S_0)$ is an element of $\basic {S_0} {n}$;
    \item given an $\mathbf S$-algebra $\mathbf T = (T_1 \stackrel{\Psi}{\to} T_0)$ and an element $y=(y_1, \dots , y_n) \in T_1^n$, the \emph{evaluation} $f(y)$ of $f$ at $y$ is the element of $T_1$ obtained by plugging in $\Psi(d^{\mathbf j} (y_i))$ for $x_i^{(\mathbf{j})}$ in $f$;
    \item $y \in T_1^n$ is a \emph{solution} for $f$ in $\mathbf T$ if $f(y)$ tropically vanishes in $T_1$;
    \item given a system $E$ of TDEs, we denote by $\textup{Sol}_{\mathbf T}(E)$ its set of solutions in $\mathbf T$.
\end{enumerate}
\end{defin}
As proved in \cite[Proposition 4.5.1]{gianmereta}, the functor of solutions to a system of TDEs can be corepresented by a quotient of the pair of differential polynomials over $\mathbf S$ by a bend congruence (see \cite{eqtrop} for a definition of bend congruence).

To conclude this section, we describe how to tropicalize a system of PDEs with coefficients in the valued differential ring $v \co R_m \rightarrow S$ to obtain the associated system of TDEs, and how to tropicalize its solutions. For systems of polynomial equations this is done via the valuation $v$, while here we need the more refined datum of a differential enhancement of a valuation. 

\begin{defin}
    Let $S$ be an idempotent semiring. We say that a map $ \widetilde v \co R_m \rightarrow S$ is a \emph{submultiplicative seminorm} if for all $a,b \in R_m$:
    \begin{enumerate}
        \item $\widetilde v(0) = 0_S$ and $\widetilde v(1) = 1_S$;
        \item $\widetilde v(a+b) \preceq \widetilde v(a) \oplus \widetilde v(b)$;
        \item $\widetilde v(ab) \preceq \widetilde v(a) \odot \widetilde v(b)$.
    \end{enumerate}
\end{defin}

\begin{defin}\label{def:diff-enhan}
For any $m$, given a valuation $ v\co R_m \to S_0$ to an idempotent semiring $S_0$, a \emph{differential enhancement of $v$} is a partial pair $\mathbf{S} = (S_1\to S_0)$ with $m$ differentials and a submultiplicative seminorm $\widetilde{v}\co R_m \to S_1$ such that:
\begin{enumerate}
\item $d_{S_1,i} \widetilde{v}(A) = \widetilde{v}(\frac{d}{dt_i}A)$ for all $A \in R_m$ and for all $i \in [m]$;
\item the following diagram commutes:
\begin{center}
\begin{tikzcd}
& S_1 \arrow[d]  \\
R_m \arrow[ur, "\widetilde{v}"] \arrow[r, swap, "v"] & S_0 .
\end{tikzcd}
\end{center}
\end{enumerate}

For a valuation $v$ together with a differential enhancement $(\mathbf S, \widetilde{v})$ we will use the term \emph{differentially enhanced valuation} and denote it as $\mathbf{v}=(v,\widetilde{v})\co R_m \to \mathbf{S}$.  
\end{defin}

We introduce here several examples of differential enhancements that we will use in the following.
\begin{ex}\label{ex:diff-enhancements}
\begin{enumerate}
\item Let $v_{\text{triv}} \co K \rightarrow \B$ be a trivially valued field. As seen in point (1) of Example \ref{ex:conv-vals}, the valuation obtained by extending the trivial valuation as in Remark \ref{oss:rank2valuation} is the $t$-adic valuation $w \co R_1\to \T$. It admits a differential enhancement on the pair $\mathbf T _1$ of point (1) of Example \ref{ex:pairs}:
\begin{center}
\begin{tikzcd}
& \B_1 \arrow[d] \\
R_1 \arrow[ur, "\widetilde{w}"] \arrow[r, swap, "w"] & \T 
\end{tikzcd}
\end{center}
in which the map $\widetilde w$ sends a power series over $K$ to its coefficientwise trivial valuation.  This is the differentially enhanced valuation used by Grigoriev \cite{grig} in his framework and subsequent works
\cite{aroca,FT20} for ODEs.
\item Let $v_K \co K \rightarrow \T$ be a nontrivially valued field and let $v \co R_1\rightarrow \T^{(2)}$ the extension of $v_K$ to $\pk$ as in Remark \ref{oss:rank2valuation}. The valuation $v$ admits a differential enhancement
\begin{center}
\begin{tikzcd}
& \pt_{v_K} \arrow[d] \\
R_1\arrow[ur, "\widetilde{v}"] \arrow[r, swap, "v"] & \T^{(2)} 
\end{tikzcd}
\end{center}
to the pair $\mathbf S _1$ of point (2) of Example \ref{ex:pairs}. Analogously as above, the map $\widetilde v$ is defined by coefficientwise application of $v_K$.
\item Given a trivially valued field $K$, the valuation $v \co R_m \to \textup{Conv}_m(\B)$ sending a multivariate power series to the convex hull of its support, introduced in point (1) of Example \ref{ex:conv-vals}, admits a differential enhancement to the pair $\mathbf T _m$ of point (3) of Example \ref{ex:pairs}:
	\begin{center}
		\begin{tikzcd}
		&\B_m \arrow[d, "\Psi_m"] \\
		R_m \arrow[ur, "\widetilde{w}"] \arrow[r,swap, "w"] &  \textup{Conv}_m(\B)
		\end{tikzcd}
	\end{center}
where $\widetilde w$ is the map $R_m \to \B_m$ given by coefficientwise trivial valuation, i.e.\ taking the support. This is the differentially enhanced valuation used in \cite{sebastian,  cotterill, sebastianmereta} in the context of PDEs. This recovers the differential enhancement of point (1) for $m=1$. We will denote the differentially enhanced valuation $(w, \widetilde w)$ as $\mathbf w$.
\item Given a valued field $v_K \co K \rightarrow \T$, consider the valuation $v \co R_m \to \textup{Conv}_m(\T)$, as defined in point (2) of Example \ref{ex:conv-vals}. This admits a differential enhancement to the pair $\mathbf S _m$ of Example \ref{ex:pairs} part (4):
	\begin{center}
		\begin{tikzcd}
		& \T_{m, v_K} \arrow[d, "\Phi_m"] \\
		R_m \arrow[ur, "\widetilde{v}"] \arrow[r,swap,"v"] & \textup{Conv}_m(\T) 
		\end{tikzcd}
	\end{center}
	where $\widetilde v$ is given by coefficientwise application of $v_K$. This recovers the differential enhancement of point (2) for $m=1$. We will denote the differentially enhanced valuation $(v, \widetilde v)$ as $\mathbf v $.
\end{enumerate}
\end{ex}
\begin{oss}\label{rem:enhance-extension}
    Let $(L,v_L)$ be a valued field extension of the valued field $(K,v_K)$. The differential enhancement $\mathbf{v}_L$ as described in point (4) above is then an extension of the differential enhancement $\mathbf{v}_K$.
\end{oss}
For any $m$, the differential enhancement of point (4) (and that of point (2)) of the example above allows us to keep track of the valuation of the coefficients of a power series in $R_m$. Notice that there is a morphism $\sigma = (\sigma_0,\sigma_1)$ of pairs between $\mathbf T _m$ and $\mathbf S _m$ of Example \ref{ex:pairs}:
\begin{equation} \label{diagram:refined-to-grig}
		\begin{tikzcd}
		\T_{m,v_K}  \arrow[d] \arrow[r, "\sigma_1"] & \B_m \arrow[d]\\
		 \textup{Conv}_m(\T) \arrow[r, "\sigma_0"] &  \textup{Conv}_m(\B)
		\end{tikzcd}
\end{equation}
where the $\sigma_0$ is the homomorphism as in (\ref{eq:sigma-0}) and $\sigma_1$ is defined by sending a tropical multivariate power series to its support. This  makes the pair $\mathbf T _m$ into an $\mathbf S _m$-algebra. 
Furthermore, $\sigma$ sends the differentially enhanced valuation $\mathbf{v}$ to the $\mathbf{w}$ of Example \ref{ex:diff-enhancements}. Thus it is the link between the framework for TDEs as developed by Grigoriev in \cite{grig} and the more general one of \cite{gianmereta}: the new framework refines Grigoriev's one and includes it as a special case, thus solutions for a differential equations in the refined setting are also solutions in Grigoriev's sense. 

\begin{oss}\label{oss:restr-of-diff-enh}
Let $\Gamma_K$ be the value group of $v_K$. As in the following we will need the differentially enhanced valuation $\mathbf{v}$ to be surjective (i.e.\ both $v$ and $\widetilde v$ to be surjective), we will consider its restriction to the pair $\mathbf{\Gamma}^K_m := (\widetilde{\Gamma}_K[\![t_1, \dots , t_m]\!] \stackrel{\Phi_m}{\to} \textup{Conv}_m(\widetilde{\Gamma}_K))$, obtained by restricting the coefficients and the weights to the value group $\Gamma_K$. The tropical pair $\mathbf{\Gamma}^K_m$ is the quotient map discussed in Remark \ref{oss:conv-as-quotient}. The differentially enhanced valuation $\mathbf{v}: R^K_m \rightarrow \mathbf{\Gamma}^K_m$ is surjective.
\end{oss}

We are finally ready to define the differential tropicalization of a system of PDEs with coefficients in $R_m$, and of the solutions of this system. Let us fix a differentially enhanced valuation 
\[
    \mathbf v =(v,\widetilde{v})\co (R_m, \left \{d/dt_i \right \}_{i \in [m]}) \to \mathbf{S}
\]
to a pair $\mathbf S$. With this data:
\begin{enumerate}
\item we tropicalize points $A \in R_m ^n$ via the map $\trop_{\widetilde v}\co R_m^n \to S_1^n$ defined by applying
$\widetilde{v}$ componentwise.
\item we tropicalize differential equations by applying $v$ coefficientwise to define a map 
\[
\trop_v\co \diff {R_m} n \to \basic {S_0} n.
\]
\end{enumerate}

In this context, given a differential ideal $I \subset \diff {R_m} n$, the statement of the fundamental theorem is the following equality of subsets of $S_1$: 
\[
	\textup{Sol}_{\mathbf S}(\trop_{v} (I)) = \trop_{\widetilde v}(\textup{Sol}_{R_m}(I)). 
\]
In case the differentially enhanced valuation considered is $\mathbf w$ of point (3) of Example \ref{ex:diff-enhancements}, i.e.\ that of Grigoriev's setting, this is the main result of \cite{aroca} for $m=1$, and of \cite{sebastian}, for any $m$. In general, the inclusion $\textup{Sol}_{\mathbf S}(\trop_{v} (I)) \subset \trop_{\tilde v}(\textup{Sol}_{R}(I))$ holds, as proved in \cite[Proposition 5.2.2]{gianmereta} (see also \cite[Proposition 7.1.5]{tesimereta} for a proof of this result in a more general formulation).

The main aim of the present work is to prove a fundamental theorem for the differential enhancement $\mathbf v$ of point (4) of Example \ref{ex:diff-enhancements}, and to develop a Gröbner theory that allows us to extend and incorporate in the present theory the results of \cite{FT20} and \cite{hugao}.

\begin{oss}
    Notice that, using the notation of Example \ref{ex:diff-enhancements} and Diagram \ref{diagram:refined-to-grig},  the following inclusion holds for every $m$:
    \[
           \sigma_1 \left ( \textup{Sol}_{\mathbf S _m}(\trop_{v} (I))  \right)
           \subset 
            \textup{Sol}_{\mathbf T _m}(\trop_{w} (I))  
    \]
    and in general it is strict, as showed in \cite[Example 4.2.2]{gianmereta}. This is evidence that the general framework developed in \cite{gianmereta} and \cite{tesimereta} refines the previous one of Grigoriev.
\end{oss}

\section{Limits of systems of algebraic equations}\label{section:preliminaries}

The main aim of this section is to prove Corollary \ref{corolla:remark+proposition} which relates solvability of systems of differential equations to solvability of certain infinite systems of algebraic equations. More precisely, given a system $F$ of differential equations over the valued field $(K,v_K)$, we produce an infinite system of polynomial equations $G$ such that $F$ has a solution in power series with coefficients in some extension of $K$ with a given valuation $S \in (\mathbb{T}_{m,v_K})^n$ if and only if every finite subset of $G$ is solvable in $K$ by a point whose valuation is equal to a certain finite tuple of coefficients of $S$. This will be essential to the proof of the Fundamental Theorem. 

The section is divided into two parts: in the first one we prove Proposition \ref{prop:compact-in-extension}, a result on infinite intersections of projections of fibers of the tropicalization. For the proof we use some tools from model theory, particularly from Hrushovski and Loeser's interpretation of Berkovich spaces as type spaces. In the second part we describe a general set-up to extract an infinite system of algebraic equations from a differential ideal, we apply Proposition \ref{prop:compact-in-extension} to show  Proposition \ref{Prop:limit-empty-iff-some-empty}, which is an analogue of \cite[Theorem 1]{bouliermodel}, and we deduce Corollary \ref{corolla:remark+proposition}. 

\subsection{Model theory and projections of fibers of the tropicalization}

 We refer to \cite[Chapter 1]{Hodges_1993} for the definitions of the basic objects of model theory, such as language, formula, structure, and of the notion of a structure $\mathcal{S}$ satisfying a formula $\phi$, which we denote (as usual) as $\mathcal{S} \models \phi$. We provide some examples which are relevant to the material in this section.

\begin{ex}\mbox{}\label{exa:languages}
    \begin{enumerate}
        \item The \emph{language of rings} $\mathcal{L}_{\textnormal{ring}}=\{0,1,+,-,\cdot\}$ contains constants $0$ and $1$ and binary functions $+,-,\cdot$. Any ring is a structure in this language by interpreting the constants and functions in the obvious way. Examples of formulas in this language are $p(\overline x)=0$ and $q(\overline{x}) \neq 0$, where $p$ and $q$ are polynomials with integer coefficients.
        \item The \emph{language of valued fields} is the language $\mathcal{L}_{\Gamma}$ with two \emph{sorts} $\textnormal{VF}$ and $\Gamma$. This means that each variable comes with an understanding of whether it stands for an element of the valued field or of the value group. The sort $\textnormal{VF}$ is equipped with the language $\mathcal{L}_{\textnormal{ring}}$, while the sort $\Gamma$ is equipped with the language $\{+,-,\leq,0 \}$ of ordered abelian groups and an additional constant symbol $\infty$ (which we interpret as the valuation of $0$), and there is a function symbol $v:\textnormal{VF} \to \Gamma$. 
    \end{enumerate}
\end{ex}

There are several possible choices for languages of valued fields. Another popular choice is a one-sorted language obtained by adding to $\mathcal{L}_{\textnormal{ring}}$ a unary relation symbol for the valuation ring.

\begin{defin}
    Let $\mathcal{S}$ be a structure on a set $S$ in a language $\mathcal{L}$. A subset $X \subseteq S^n$ is \emph{definable} if there is an $\mathcal{L}$-formula $\phi(\bar x)$ such that $X=\{a \in S^n \mid \mathcal{S} \models \phi(a)\}.$
\end{defin}

\begin{ex}
    \mbox{}\label{exa:definable-sets}
    \begin{enumerate}
        \item If $\mathcal{S}$ is the complex field as a structure in $\mathcal{L}_{\textnormal{ring}}$, then by \emph{quantifier elimination} for algebraically closed fields \cite[Theorem 2.7.3]{Hodges_1993} every formula is equivalent to a quantifier-free formula, namely to a Boolean combination of polynomial equations and inequations. Hence a subset of $\mathbb{C}^n$ is definable if and only if it is constructible.
        \item Examples of definable sets in an algebraically closed valued field, seen as a structure in $\mathcal{L}_{\Gamma}$, include every set of the form $\{x \in V \mid v(x)=\gamma \}$ where $V$ is an algebraic subvariety of $K^n$ and $\gamma$ an $n$-tuple in the value group.  
    \end{enumerate}
\end{ex}

Definable sets can be seen as functors from the category of models of a given theory (with elementary embeddings as maps) to the category of sets. Given an extension $(K,v_K) \leq (L,v_L)$ of algebraically closed valued fields and a definable subset $X \subseteq K^n$ defined by some formula $\phi$, we write $X(L)$ for the set of elements of $L^n$ which satisfy $\phi$. 

\begin{defin}
    Let $\mathcal{S}$ be a structure on a set $S$ in a language $\mathcal{L}$, $A\subseteq S$ a subset. A \emph{type}  over $A$ in $n$ variables is a set $p$ of formulas in $n$ variables with parameters from $A$ such that:
    \begin{itemize}
        \item For every $\phi_1,\dots,\phi_k \in p$, we have $\mathcal{S} \models \exists x_1,\dots,x_n \bigwedge_{i=1}^k \phi_i(x_1,\dots,x_n)$.
        \item For every formula $\phi$ in $n$ variables in the language of $\mathcal{S}$, one of $\phi$ and $\neg \phi$ lies in $p$.
    \end{itemize}

    If there is a tuple $a \in S^n$ such that $\mathcal{S} \models \phi(a)$ for every formula $\phi \in p$, then we say $p$ is \emph{realized} in $\mathcal{S}$ by $a$.
\end{defin}

\begin{ex}\label{exa:type-examples}
   We provide some examples in the languages of Example \ref{exa:languages}:
   \begin{enumerate}
       \item If $\mathcal{S}$ is the complex field, as a structure in the language of rings, and $A=\varnothing$, then every monic polynomial $p(x)$ in one variable over $\mathbb{Q}$ which is irreducible over $\mathbb{Q}$ determines a type, the type of a root of $p(x)$ in $\mathbb{C}$. The fact that any two roots of the same polynomial (over $\mathbb{Q}$ and irreducible over $\mathbb{Q}$) satisfy the same type follows from the fact that there is an automorphism of $\mathbb{C}$ which maps one to the other.
       \item If $\mathcal{S}$ is a valued field, seen as a structure in the language $\mathcal{L}_{\Gamma}$, we may consider the \emph{generic type} $p_R$ of the valuation ring. This is the type which contains all formulas of the form $v(f(x))=0$ where $f$ is a polynomial over the valuation ring of $K$. Since the Gauss extension of any valued field is unique up to isomorphism, this determines a complete type over $K$. Clearly $p_R$ is not realized in $K$, since for every element $a$ in the valuation ring of $K$ the formula $v(x-a)=0$ belongs to $p_R$ and therefore $a$ does not realize $p_R$.  \label{exa:generic-type}
   \end{enumerate}
\end{ex}

    We now move to considerably more advanced material and recall some notions and results around the interplay between Berkovich spaces and the model theory of algebraically closed valued fields \cite[Chapter 14]{HrushovskiLoeser}.
    
    \begin{defin}[{\cite[p.18]{HrushovskiLoeser}
    }]
        Let $(K,v_K)$ be an algebraically closed valued field, embedded in an algebraically closed valued field extension with valuation surjective on $\mathbb{T}$. A type $p$ over $K \cup \mathbb{T}$ is \emph{almost orthogonal} to $\Gamma$ if there is an algebraically closed valued field extension $(L,v_L)$ of $(K,v_K)$ with value group contained in $\mathbb{T}$ in which $a$ is realized.
    \end{defin}

    \begin{ex}
        As an example, we may consider the generic type of the valuation ring introduced in Example \ref{exa:type-examples}(\ref{exa:generic-type}). Let $p_R$ be the generic type of the valuation ring of some algebraically closed valued field $K$, and let $a$ be a realization of $p_R$ in an algebraically closed valued field extension $L$ of $K$. We may assume that $L$ is the algebraic closure of $K(a)$. Then the extension $L/K$ has transcendence degree 1, and so does the residue field extension by definition of the generic type. By the Zariski-Abhyankar inequality \cite[Corollary 3.25]{vandenDries2014} the image of the valuation on $L$ coincides with the image of the valuation on $K$, so in particular it is contained in $\mathbb{T}$.

        More generally, whenever $(K,v_K) \subseteq (L,v_L)$ is an extension of algebraically closed valued fields, the value group of $L$ is contained in $\mathbb{T}$ (hence so is the value group of $K$) and $a \in L^n$, then set of formulas $\phi$ with parameters from $K$ such that $L \models \phi(a)$ forms a type over $K$, almost orthogonal to $\Gamma$. This is the case, for example, if $K$ is the algebraic closure of $\mathbb{Q}_p$ and $L$ is $\mathbb{C}_p$.

        On the other hand, if a type $p$ contains the set of formulas $\{x>n\mid n \in \N\}$ then the valuation of any realization of $p$ is larger than every real, so $p$ is not almost orthogonal to $\Gamma$.
    \end{ex}

\begin{defin}
    Let $(K,v_K)$ be an algebraically closed valued field, $X \subseteq K^n$ a definable set. We write $B_K(X)$ for the set of types over $K \cup \mathbb{T}$ which are almost orthogonal to $\Gamma$ and contain the formula expressing $x \in X$. 
    
    We endow $B_K(X)$ with the topology given by taking as the pre-basic open sets the subsets $\{p \in B_K(X) \mid ``x \in U" \in p\}$ and $\{p \in B_K(X) \mid ``v_K(f(x)) \in W" \in p\}$ as $U,f,W$ vary among Zariski open subsets of $X$, regular functions on $U$, and open subsets of $\mathbb{T}$ respectively.
\end{defin}

As observed in \cite[Chapter 14]{HrushovskiLoeser}, when $X$ is an algebraic variety, $B_K(X)$ can be identified with the underlying topological space of the Berkovich analytification of $X$. Note that $B_K(X)$ is Hausdorff for every definable $X \subseteq K^n$.

\begin{prop}\label{prop:compact-in-extension}
Let $\{X_i \mid i \in \N\}$ be a decreasing family of nonempty subsets of $K^n$ such that for every $i$ there are $\ell_i \in \N$, an algebraic variety $V_i \subseteq K^{n+\ell_i}$, and $\gamma_i \in \mathbb{T}^{n+\ell_i}$ such that \[X_i=\{x\in K^n \mid \exists y \in K^{\ell_i} \text{ s.t. } (x,y) \in V_i \text{ and } v_K(x,y)=\gamma_i \}.\] Then there is an extension $(L,v_L)$ of $(K,v_K)$, with value group contained in $\mathbb{T}$, such that the intersection $\bigcap_{i \in \N} X_i(L)$ is non-empty. 
\end{prop}

\begin{proof}
    For every $i \in \N$, the set $V_i \cap v_K^{-1}(\gamma_i)$ is \emph{bounded} (it is contained in a valuative ball centered at the origin) and it is closed for the \emph{$v+g$-topology} introduced in \cite[Section 3.7]{HrushovskiLoeser}: this means that it is closed both for the valuation topology and for the \emph{$g$-topology} \cite[Definition 3.7.1]{HrushovskiLoeser}. To see it is closed for the $g$-topology, we observe that its complement is defined by the formula \[x \notin V_i \vee \left(x \in V_i \wedge \bigvee_{j=1}^{n+\ell_i} (v_K(x_j) >(\gamma_i)_j \vee v_K(x_j) < (\gamma_i)_j)\right)\] namely a positive Boolean combination of Zariski open and closed sets and sets defined by valuative inequalities, which is precisely the definition of a $g$-open set.

    It then follows by \cite[Proposition 14.1.2]{HrushovskiLoeser} that $B_K(V_i \cap v_K^{-1}(\gamma_i))$ is compact for every $i$. Consider the natural projection map $\pi_i:B_K(V_i \cap v_K^{-1}(\gamma_i)) \to B_K(X_i)$ which maps the type $p \in B_K(V_i)$ to the set of the formulas in $p$ which only involve the first $n$ variables. This is a surjective function: given a type $p$ which is almost orthogonal to $\Gamma$ and contains the formula $x \in X_i$, it has a realization $a \in L^n$ for some extension $(L,v_L)$ of $(K,v_K)$. By definition of $X_i$, there is $b \in L^{\ell_i}$ such that $(a,b) \in V_i$ and $v_K(a,b)=\gamma_i$. The type $q$ over $K \cup \mathbb{T}$ of $(a,b)$ (that is, the set of formulas $\phi$ in $n+\ell_i$ variables with parameters in $K \cup \mathbb{T}$ such that $L \models \phi(a,b)$) is then almost orthogonal to $\Gamma$, so it belongs to $B_K(V_i\cap v_K^{-1}(\gamma_i))$, and $\pi_i(q)=p$. 
    
    It easily follows from the definition of the topology on $B_K(V_i \cap v_K^{-1}(\gamma_i))$ that $\pi_i$ is continuous, so the image of $\pi_i$ is compact; hence, every $B_K(X_i)$ is compact. Since $B_K(X_0)$ is Hausdorff, each $B_K(X_i)$ is a closed subset of $B_K(X_0)$, so $\bigcap_{i \in \N} B_K(X_i)$ is a decreasing intersection of closed subsets of a compact space which has the finite intersection property, so it is nonempty. By definition of almost orthogonality to $\Gamma$, any type in this intersection is realized in an extension $(L,v_L)$ of $(K,v_K)$ with image of the valuation contained in $\T$ by some element $\alpha \in L^n$. Thus, $\alpha \in X_i(L)$ for every $i \in \N$.
    \end{proof}

    \begin{oss}
        For the model-theoretically inclined reader, we stress that while a version of Proposition \ref{prop:compact-in-extension} in which we allow for valuations of arbitrary rank is of course a trivial consequence of the compactness theorem, in this paper we only consider valuations of rank $1$. Realizing the intersection of the $X_i$ by compactness does not guarantee any control on the value group of the extension. We work with types almost orthogonal to $\Gamma$ to make sure that the partial type defining the intersection has a completion whose realizations do not add new archimedean classes to the value group.
    \end{oss}

\subsection{Solving infinite systems in extensions}
 
Having established Proposition \ref{prop:compact-in-extension}, we now aim to apply it to the projections of fibers of tropicalizations of a differential ideal's truncations.
Let us fix here the notation for the following. 

Fix $m >0$ an integer. For $\jj \in \N^m$ let $\jj ! := \prod_{i =1}^m j_i!$ and $\norm \jj = \sum_{i=1}^m j_i$.  Let $v \co R_m \rightarrow \textup{Conv}_m(\T)$ be the valuation constructed from $v_K$ as in Example \ref{ex:conv-vals}. Fix the differentially enhanced valuation $\mathbf{v}=(v,\widetilde{v})\co (R_m,\{d/dt_i\}_{i \in [m]}) \to \mathbf{S_m}$ of Example \ref{ex:diff-enhancements}. 

From the trivial valuation $v_{\text{triv}} \co K \rightarrow \B \ \subset \T$ on $K$ instead we construct the valuation $w \co R_m \rightarrow \textup{Conv}_m(\B)$. Fix the differentially enhanced valuation \[
    \mathbf{w}=(w,\widetilde{w})\co (R_m,\{d/dt_i\}_{i \in [m]}) \to \mathbf{T}_m
\]
of Example \ref{ex:diff-enhancements}: this is the differential enhancement considered in Grigoriev's framework. Via the morphism $\sigma$ of diagram \ref{diagram:refined-to-grig} we endow $\mathbf S _m$ with the structure of a $\mathbf T _m$-algebra and we relate the two differential enhancements $\mathbf{w}$ and $\mathbf{v}$. We start by recalling some notation and objects from \cite{aroca} and \cite{sebastian} and extend their definition to the nontrivally valued case.

Let $I \subset \diff {R_m} n$ be a differential ideal. By \cite[page 21]{ritt}, there exist finitely many differential polynomials $f_1, \dots , f_s \in I$ such that the set of solutions to $I$ is equal to the set of solutions of $f_1, \dots , f_s$.

For $1 \le l \le s$ and $\rr \in \N^m$, set:
\[
        F_{l,\mathbf{r}} :=(d^\rr f_l) \rvert _{t=0} \in 
        K[x_i^{(\jj)} \mid i=1, \dots, n ; \mathbf{j} \in \N^m]
\]
and 
\[
    A_\infty := V \left (
    \{F_{l,\mathbf{r}}\}_{\substack{1\le l \le s\\ \mathbf{r} \in \N^m}} \right ) \subset \left ( K^{\N^m} \right )^n
\]
The map $\Xi \co \left ( K^{\N^m} \right )^n \rightarrow R_m^n$ defined as
\[
    a := \left ( (a_{1,\jj})_{\jj \in \N^m}, \dots , (a_{n,\jj})_{\jj \in \N^m} \right ) \mapsto
     \left ( \sum_{\jj \in \N^m} \frac{1}{\jj!} a_{1,\jj} t^\jj, \dots ,\sum_{\jj \in \N^m} \frac{1}{\jj!} a_{n,\jj} t^\jj\right )
\]
is a bijection. Furthermore, as proved in \cite[Lemma 6.2]{sebastian}, given a differential polynomial $f \in \diff {R_m} n$ and $a \in \left ( K^{\N^m} \right )^n $ the following equality holds:
\[
        f \left (\Xi(a) \right ) = \sum_{\rr \in \N^m} \left ( \frac{1}{\rr!}\left (d^\rr(f) \right ) \rvert _{t=0}(a) \right) t^\rr .
\]
Thus we obtain:
\[\textup{Sol}_{R_m}(I) = \Xi(A_\infty).\]

For $k \in \N$, let $J_k$ be the smallest natural number such that 
\[
         F_{l,\rr} \in 
        K[x_i^{(\jj)} \mid i=1, \dots, n ; \norm \jj \le J_k] \quad \text{ for all } 1 \le l \le s , \,  \norm \rr \le k
\]
The cardinality of the set $\{\jj \in \N^m \mid \norm \jj \le J_k\}$ is $\binom{J_{k}+m}{m}$. Let us denote this number as $N_k$ and let
\[
        A_k:= V \left ( \{F_{l,\rr}\}_{\substack{1\le l \le s\\ \norm \rr \le k} } \right ) \subset \left ( K^{N_k} \right )^n.
\]
For $k \ge k' \ge 0$, let $\pi_{(k,k')} \co \left ( K^{N_k} \right )^n \rightarrow \left ( K^{N_{k'}} \right )^n$ be the projection morphism forgetting the last $N_k - N_{k'}$ entries of every vector.
We have
\[
    \pi_{(k,k')}(A_k) \subset A_{k'}
\]
and $A_\infty$ is the inverse limit of the system given by the sets $A_k$ and the maps $\pi_{(k,k')}$:
\[
    A_\infty = \varprojlim A_k.
\]

	\begin{defin}
        Let us denote again by $\sigma_1$ the restriction $\T \rightarrow \B$ of $\sigma_1$ to $\T$.
        Let  $k \in \N$ and $S := (S_1, \dots , S_n) \in (\T_m)^n$, where we write $S_i$ as $\sum_{\jj \in \N^m}{c_{i, \jj}} t^\jj $ for every $i = 1, \dots , n$. With this notation, define: 
		\[
			(\mathbb{V}_\infty)_S^{v_{\triv}} := \left  \{(x_{i,\jj})_{\substack{i = 1, \dots, n \\ \jj \in \N^m}} \in \left ( K^{\N^m } \right )^n \bigmid v_{\triv}(x_{i,\jj}) = \sigma_1(c_{i,\jj}) \text{ for all } i,\jj \right   \}
		\]
		\[
			(\mathbb{V}_\infty)_S^{v_K}:= \left  \{
            (x_{i,\jj})_{\substack{i = 1, \dots, n \\ \jj \in \N^m}} \in \left ( K^{\N^m} \right )^n  \bigmid  v_K(x_{i,\jj}) = c_{i,\jj} + v_K(\jj!) \text{ for all } i,\jj \right \}
		\]
  and 
        \[
			(\mathbb{V}_k)_S^{v_{\triv}} := \left  \{(x_{i,\jj})_{\substack{i = 1, \dots, n \\ \norm \jj \le J_k}} \in \left ( K^{N_k} \right )^n \bigmid v_{\triv}(x_{i,\jj}) = \sigma_1(c_{i,\jj}) \text{ for all } i,\jj \right   \}
		\]
		\[
			(\mathbb{V}_k)_S^{v_K} := \left  \{(x_{i,\jj})_{\substack{i = 1, \dots, n \\ \norm \jj \le J_k}} \in \left ( K^{N_k } \right )^n  \bigmid  v_K(x_{i,\jj}) = c_{i,\jj} + v_K(\jj!) \text{ for all } i,\jj \right \}
         \]
    Furthermore, let 
    \[
    \left ( A_\infty \right )^{v_\triv}_S := A_\infty \cap (\mathbb{V}_\infty)_S^{v_{\triv}}
    \quad \quad 
    \left ( A_\infty \right )^{v_K}_S := A_\infty \cap (\mathbb{V}_\infty)_S^{v_K}
    \]
    and for every $k \in \N$
    \[
    \left ( A_k \right )^{v_\triv}_S := A_k \cap (\mathbb{V}_k)_S^{v_\triv}  
    \quad \quad 
    \left ( A_k \right )^{v_K}_S := A_k \cap (\mathbb{V}_k)_S^{v_K}. 
    \]
	\end{defin}
 
	\begin{oss} \label{remark:trop-fibers-contained}
    The sets $(\mathbb{V}_k)_S^{v_{\triv}}$ and $(\mathbb{V}_k)_S^{v_K}$ are the fibers of the tropicalization with respect to $v_\text{triv}$ and $v_K$, respectively.
	Furthermore, for every $k \in \N$ the following inclusion holds:
	\[
			(\mathbb{V}_k)_S^{v_K}  \subseteq (\mathbb{V}_k)_S^{v_\text{triv}}  
	\] 
	and $(\mathbb{V}_k)_S^{v_\text{triv}} $ is a torus contained in $\left ( K^{N_k} \right )^n$, of dimension equal to 
	\[|\{ (i,\jj) \mid c_{i,\jj} \neq \infty \text{ and } \norm \jj \le J_k\}|.\]
 \end{oss}

\begin{oss} \label{oss:solutions-A}
 In the case of trivial valuation, an element $S=(S_1, \dots , S_n) \in (\B_m)^n$ is in  $\trop_{\widetilde w}(\textup{Sol}_{R_m}(I))$ if and only if there exists $a \in A_\infty$ such that $\trop_{\widetilde w}(\Xi(a)) = S$, i.e.\ if the set $\left ( A_\infty \right )^{v_\triv}_S $ is non-empty. 
 Analogously, in the nontrivially valued case, $S=(S_1, \dots , S_n) \in (\T_{m,v_K})^n$ is in  $\trop_{\widetilde v}(\textup{Sol}_{R^K_m}(I))$ if and only if the set $\left ( A_\infty \right )^{v_K}_S $ is non-empty. 
\end{oss}

For $k \ge k' \ge 0$, the following inclusions hold:
\[
    \pi_{(k,k')} \left ( (\mathbb{V}_k)_S^{v_{\triv}} \right ) \subseteq (\mathbb{V}_{k'})_S^{v_{\triv}} 
    \quad \quad 
    \pi_{(k,k')} \left ( (\mathbb{V}_k)_S^{v_K} \right ) \subseteq (\mathbb{V}_{k'})_S^{v_K} 
\]
and 
\[
    \pi_{(k,k')} \left ( \left ( A_k \right )^{v_\triv}_S \right ) \subseteq \left ( A_{k'} \right )^{v_\triv}_S
    \quad \quad 
    \pi_{(k,k')} \left ( \left ( A_k \right )^{v_K}_S \right ) \subseteq \left ( A_{k'} \right )^{v_K}_S
\]
The set $(A_\infty)_S^{v_K}$ (respectively $(A_\infty)_S^{v_\triv}$) is the inverse limit of the system given by the sets $(A_k)_S^{v_K}$ (respectively $(A_k)_S^{v_\triv}$) and the maps $\pi_{(k,k')}$. 

\begin{lem}\label{lem:binomials-then-subgroups}
    Let $I \subseteq R^K_m\{x_1,\dots,x_n\}$ be a differential ideal, and assume $I$ is generated as a differential ideal by a finite set of linear differential binomials with constant coefficients. Then for every $k \in \N$ the set $A_k \cap ((K^{\times})^{N_k})^n$ is a coset of an algebraic subgroup of $((K^{\times})^{N_k})^n$.
\end{lem}
\begin{proof}
   Assume $I$ is generated as a differential ideal by a finite set of linear binomials $x_{i_1}^{(\jj_1)}-\lambda x_{i_2}^{(\jj_2)}$ where $\lambda \in K$ is a constant coefficient (depending on $i_1,\jj_1,i_2,\jj_2$). For every $\mathbf{r} \in \N^m$, the derivative $d^{\mathbf{r}} (x_{i_1}^{(\jj_1)}-\lambda x_{i_2}^{(\jj_2)})$ equals $x_{i_1}^{(\jj_1+\mathbf{r})}-\lambda x_{i_2}^{(\jj_2+\mathbf{r})}$, so each of the polynomials $F_{\mathbf{r}}$ defined as above is a binomial. Hence $A_k$ is defined by binomials, so its intersection with $((K^{\times})^{N_k})^n$ is a coset of an algebraic subgroup of $((K^{\times})^{N_k})^n$.
\end{proof}

\begin{prop} \label{Prop:limit-empty-iff-some-empty}
		Let $I \subseteq R^K_m\{x_1,\dots,x_n\}$ be a differential ideal, and let $S \in (\Gamma_K)_{m}^n$. 
        \begin{enumerate}
            \item $(A_k)^{v_K}_S$ is nonempty for all $k \in \N$ if and only if there is an extension $(L,v_L)$ of $(K,v_K)$, with a differential enhancement of $v_L$, such that $(A_\infty)^{v_L}_S$ is nonempty.
            \item  Assume $I$ has a set of generators consisting of linear differential binomials with constant coefficients. Then  $(A_k)^{v_K}_S$ is nonempty for all $k \in \N$ if and only if $(A_{\infty})^{v_K}_S$ is nonempty.
        \end{enumerate}
	\end{prop}

    \begin{proof}
    We first prove Point (1).
    
        $(\Leftarrow)$ Assume there exists such an extension $(L,v_L)$, with a differential enhancement of $v_L$, and let $a \in (A_\infty)^{v_L}_S$. Then for every $k \in \N$ we have $\pi_k(a) \in (A_k)_S^{v_L}$, so $v_L(\pi_k(a))$ lies in the tropicalization of $A_k(L)$; observe that $v_L(\pi_k(a)) \in (\Gamma_K)^{N_k}$, since the coefficients of $S$ lie in $\Gamma_K$. Since $A_k$ is defined over $K$, by \cite[Lemma 2.6.5 and Theorem 2.6.6]{macsturm} there is a finite set $\mathcal{G}$ of polynomials with coefficients in $K$ such that the tropicalization of $A_k(L)$ is the intersection of the tropical hypersurfaces of the polynomials in $\mathcal{G}$. This implies that $\trop(A_k)=\trop(A_k(L))$, so that $v_L(\pi_k(a)) \in \trop(A_k)$; by the Fundamental Theorem of tropical geometry \cite[Theorem 3.2.3]{macsturm} we then have that $(A_k)_S^{v_K} \neq \varnothing$.

        $(\Rightarrow)$ Assume that $(A_k)_S^{v_K} \neq \varnothing$ for every $k$. By Proposition \ref{prop:compact-in-extension}, there is an extension $L$ in which $\bigcap_{i \in \N} \pi_{(i,0)}((A_i)_S^{v_L})(L) \neq \varnothing$. The valuation on $L$ admits a differential enhancement extending the differential enhancement of $v_K$, as observed in Remark \ref{rem:enhance-extension}. By \cite[page 198, Proposition 5]{bourbaki}, it follows that $(A_{\infty})_S^{v_L} \neq \varnothing$. 

        For Point (2), we only prove the right-to-left direction, the other one being obvious. If the generators of $I$ are linear differential binomials then by Lemma \ref{lem:binomials-then-subgroups} the sets $A_k \cap ((K^\times)^{N_{k}})^n$ are cosets of algebraic subgroups of $((K^\times)^{N_{k}})^n$, so, as $v_K$ is a group homomorphism, for every $k \in \N$ we have that $\trop(A_k) \cap (\mathbb{R}^{N_k})^n$ is an affine subspace of $(\mathbb{R}^{N_k})^n$. Since the class of cosets of algebraic subgroups is closed under projections and (infinite) intersections, we have that $\bigcap_{i \in \N} \pi_{(i,0)}(A_i)$ is a coset of an algebraic subgroup of $((K^\times)^{N_0})^n$.
        Let $(a_{i,\jj})_{i=1,\dots,n;  ||\jj|| \leq J_0} \in \pi_{(k,0)}(A_k) \cap (\mathbb{V}_0)_{S}^{v_K}$ and consider the algebraic variety \[C_k:=A_k \cap \{x_{i,\jj}=a_{i,\jj}| i=1,\dots,n; \, \norm \jj \leq J_0\}.\] It is easily proved by induction that 
        \[
            \trop(C_k)=\trop(A_k) \cap \{x_{i,\jj}=v_K(a_{i,\jj}) \mid i=1,\dots,n; \, \norm \jj  \leq J_0  \}. 
        \]
            By assumption, the tuple $(c_{i,\jj}+v_K(\jj!))_{i=1,\dots,n; \norm \jj \leq J_k}$ lies in $\trop(C_k)$, so by the Fundamental Theorem of tropical geometry \cite[Theorem 3.2.3]{macsturm} there is a tuple $(a_{i,\jj})_{i=1,\dots,n; \, J_0 < \norm \jj \leq J_k}$ which satisfies 
            \[
                v_K (a_{i,\jj})=c_{i,\jj}+v_K(\jj!) 
            \]
            for every $i=1,\dots,n$ and $\jj$ with $J_0 < \norm \jj \leq J_k$ and 
            \[
                (a_{i,\jj})_{\substack{i=1,\dots,n\\ \, \norm \jj \leq J_k}} \in A_k. 
            \]
            This proves that $\pi_{(k,0)}(A_k) \cap (\mathbb{V}_0)_S^{v_K} \subseteq \pi_{(k,0)}((A_k)_S^{v_K})$. Hence, 
            \[
                \bigcap_{i \in \N} (\pi_{(i,0)}(A_i)_S^{v_K})=\bigcap_{i \in \N}(\pi_{(i,0)}(A_i) \cap (\mathbb{V}_0)_S^{v_K})=\left(\bigcap_{i \in \N}(\pi_{(i,0)}(A_i))\right) \cap (\mathbb{V}_0)_S^{v_K}.
            \]
            We have already observed  that $\bigcap_{i \in \N}(\pi_{(i,0)}(A_i))$ is a coset of an algebraic subgroup. Moreover, by \cite[Corollary 3.2.13]{macsturm} we have that, for every $k \in \N$, the equality $\trop(\pi_{(k,0)}(A_k))=\pi_{(k,0)}(\trop(A_k))$ holds (abusing notation so that $\pi_{(k,0)}$ also denotes the projection map $(\mathbb{T}^{N_k})^n \to (\mathbb{T}^{N_0})^n$). It follows that 
            \[
                (c_{i,\jj}+v_K(\jj!))_{\substack{i=1,\dots,n \\ \norm \jj \leq J_0}} \in \trop \left(\bigcap_{i \in \N} \pi_{(i,0)}(A_i) \right)
            \]
            so we can conclude, again by the Fundamental Theorem, that \[\left(\bigcap_{i \in \N}(\pi_{(i,0)}(A_i))\right) \cap (\mathbb{V}_0)_S^{v_K} \neq \varnothing\] and applying \cite[page 198, Proposition 5]{bourbaki} we complete the proof.
    \end{proof}

    \begin{cor}
        Let $(K,v_K)$ be an algebraically closed valued field with a differentially enhanced valuation. There is an algebraically closed valued field extension $(L,v_L)$ of $(K,v_K)$, with a differentially enhanced valuation, such that for every differential ideal $I \subseteq R^L_m\{x_1,\dots,x_n\}$ and every $S \in (\Gamma_L)_m^n$, we have that $(A_k)_S^{v_L} \neq \varnothing$ for every $k \in \N$ if and only if $(A_\infty)_S^{v_L} \neq \varnothing$.
    \end{cor}
    \begin{proof}
        Fix an enumeration $\{K_i \mid i \in \kappa\}$, for an appropriate cardinal $\kappa$, the algebraically closed field extensions of $K$ with a differential enhancement of the valuation. The relation of being an extension of differentially enhanced valued fields naturally gives a partial order on this set. This partial order is closed under unions of chains, so by Zorn's lemma it has a maximal element $(L,v_L)$. 

    Now let $I \subseteq R_m^L\{x_1,\dots,x_n\}$ be a differential ideal and let $S \in (\Gamma_L)_m^n$. By  Proposition \ref{Prop:limit-empty-iff-some-empty}, there is an algebraically closed valued field $(L',v_{L'})$ with a differentially enhanced valuation such that if $I$ gives rise to sets $A_k$ such that $(A_k)_{S}^{v_L} \neq \varnothing$, then $(A_{\infty})_{S}^{v_{L'}} \neq \varnothing$. By maximality of $L$, we have $L=L'$, so $(A_{\infty})_S^{v_L} \neq \varnothing$ as required.
    \end{proof}

    \begin{cor}\label{corolla:remark+proposition}
        Let $(K,v_K)$ be an algebraically closed valued field with a differential enhancement of the valuation. Then there is an algebraically closed valued field extension $(L,v_L)$ of $(K,v_K)$ with differentially enhanced valuation such that for every differential ideal $I \subset R^L_m\{x_1,\dots,x_n\}$ and every $S \in (\T_{m,v_L})^n$, we have that  $(A_k)_S^{v_L} \neq \varnothing$ for all $k \in \mathbb{N}$ if and only if  $S$ lies in the set $\trop_{\tilde v}(\textup {Sol}_{R^L_m}(IR^L_m\{x_1, \dots x_n\}))$ .
    \end{cor}

    \begin{proof}
        It follows from Remark \ref{oss:solutions-A} and Proposition \ref{Prop:limit-empty-iff-some-empty}.
    \end{proof}

    By Proposition \ref{Prop:limit-empty-iff-some-empty}(2), similar proofs show that if we work with a single ideal $I$ which has a basis of linear differential binomials with constant coefficients then there is no need to pass to the extension $L$. However, it is necessary in general to pass to an extension, as shown by the following example.

    \begin{ex}\label{ex:Puiseux}
        Let $K:=\bar{\mathbb{Q}}\{\!\{s\}\!\}$ be and $L:=\mathbb{C}\{\!\{s\}\!\}$ be the fields of Puiseux series over the algebraic and complex numbers respectively. By \cite[Corollary 4.7]{DenefLipschitz}, there is a differential ideal $I \subseteq \bar{\mathbb{Q}}[\![t]\!]\{x\}$ which has a zero $C=\sum_{i=0}^\infty a_it^i \in \mathbb{C}[\![t]\!]$, but no solutions in $\bar{\mathbb{Q}}[\![t]\!]$. 

        The series $C$ can be seen as an element of the field $L[\![t]\!]$ and it is therefore a zero of $IL[\![t]\!]\{x\}$. Every nonzero coefficient of $C$ has valuation $0$, hence there is a boolean power series $S \in \mathbb{B}_1$ such that $S \in \trop_{\tilde v}(\textup {Sol}_{L[\![t]\!]}(I))$, so, by Corollary \ref{corolla:remark+proposition}, we have that $(A_k)_S^{v_K} \neq \varnothing$ for every $k$. However, if $S\in \trop_{\tilde v}(\textup {Sol}_{K[\![t]\!]}(I)) $, then there is a power series $\sum_{i=0}^\infty b_i(s)t^i$ that is a zero of $IK[\![t]\!]\{x\}$, whose coefficients $b_i(s)$ are Puiseux series over $\bar{\mathbb{Q}}$ with valuation in $\mathbb{B}$: that is, for every $i$, either $b_i(s)=0$ or $b_i(0) \in \bar{\mathbb{Q}}^\times$. It follows that the power series $\sum_{i=0}^{\infty} b_i(0)t^i $ is a zero of $I$ in the power series ring $\bar{\mathbb{Q}}[\![t]\!]$, which by assumption does not exist.
    \end{ex}
    
	\section{Statement and proof of the first part of the fundamental theorem}\label{section:theorem}
    
    This section is devoted to the proof of the first equality between the sets appearing in Theorem \ref{theorem:fundamental-introduction}: 
\begin{teorema} [Fundamental theorem of tropical differential algebra]\label{theorem:fundamental}
Let $(K,v_K)$ be an algebraically closed valued field of characteristic 0. There is an algebraically closed valued field $(L,v_L)$ extending $(K,v_K)$, with a surjective differentially enhanced valuation 
        \[
            \mathbf{v}=(v,\widetilde{v})\co (R_m^L , \{d/dt_i\}_{i \in [m]}) \to \mathbf{\Gamma}_m^L
        \]
as in Remark \ref{oss:restr-of-diff-enh} such that for every differential ideal $I$ in $R_m^K\{x_1,\dots,x_n\}$ the following equality holds: 
	\[
        \textup{Sol}_{\mathbf{ \Gamma}^L_m}(\trop_{v} (IR_m^L\{x_1,\dots,x_n\})) = \trop_{\widetilde v}(\textup{Sol}_{R^L_m}(IR_m^L\{x_1,\dots,x_n\})).
		\]
\end{teorema}

We introduce here the tropical counterpart to the map $\Xi$ defined previously. It will be useful in the proof of Theorem \ref{theorem:fundamental} above.
	\begin{defin}
	Let $\Xi_{\trop} \co \T^{\N^m} \rightarrow \T_{m,v_K}$ be the map defined by
	\[
	\Xi_{\trop}((b_{\jj})_{\jj \in \N^m} ) = \sum_{\jj \in \N^m}(b_\jj - v_K(\jj!)) t^\jj.
	\]
	It is bijective with inverse defined as follows: 
	\[
		\Xi_{\trop}^{-1} (S) = \left ((d^\jj_{v_K} S)|_{t=\infty} \right )_{\jj\in \N^m}.
	\]
        We denote again by $\Xi_{\trop}$ the map $(\T^{\N^m})^n \rightarrow \T_{m,v_K}^n$ obtained by applying $\Xi_{\trop}$ coordinatewise.
    \end{defin}
	\begin{oss}
	    For every valued field $(K,v_K)$, the  following diagram commutes: 
	\begin{center}
		\begin{tikzcd}
			K^{\N^m} \arrow[r, "\Xi"] \arrow[d,  swap, "v_K"]& R^K_m \arrow[d, "\tilde v"] \\
			\T^{\N^m} \arrow[r, swap, "\Xi_\trop"] & \T_{m,v_K}.
		\end{tikzcd}
	\end{center}
	\end{oss}
Let us set some notation for the following. 
For a differential polynomial $f \in \diff {R_m}{n}$, the \emph{order} of $f$ is the maximum of $\norm \jj$ such that there exists $i \in \{1, \dots , n\}$ such that $x_i^{(\jj)}$ is in the support of $f$.

Let $f \in \diff {R_m}{n}$ be of order less than or equal to $r \in \N$, and let us write it as $f(x) = \sum_{\lambda \in \Lambda} A_\lambda x^\lambda$ where 
$\Lambda \subset\text{Mat}_{\binom{r+m}{m} \times n}(\N)$
is finite, $A_\lambda$ is an element of $R_m$ for every $\lambda$ and $x^\lambda$ is the differential monomial defined as 
\[  
    \prod_{\substack{i \in \{1, \dots, n\} \\ \norm \jj \le r}}  ( x_i^{(\jj)}  ) ^{\lambda_{i,\jj}}. 
\]
We will denote by $C_\lambda$ the element $v(A_\lambda) \in \textup{Conv}_m(\T)$. Thus, the tropicalization of $f$ is 
\[
        \sum_{\lambda \in \Lambda} C_{\lambda} x^\lambda
        \in \basic {\textup{Conv}_m(\T)} n
\] 
for the same $\Lambda$. As in Remark \ref{oss:representation-of-conv}, let $\ell_\lambda$ be the number of vertices of $C_\lambda$. We can represent $C_\lambda$ as 
\[
    C_\lambda = ((\jj_{\lambda,1};\alpha_{\lambda, 1}), \dots ,(\jj_{\lambda,\ell_{\lambda}};\alpha_{\lambda, \ell_{\lambda}}))
\]
where the elements of the set $\{\jj_{\lambda,i}\}_{i \in [\ell_\lambda]}$ are in the support of $A_{\lambda}$ and $\alpha_{\lambda, i} \in \R$ for all $i \in [\ell_\lambda]$.

Given a differential monomial $M_\lambda$ of order $\le r$, we will denote as $M_\lambda(S)$ its evaluation at an element $S \in \T_{m,v_K}$ as in Definition \ref{defin:tropical-solutions}, when we see $M_\lambda$ as an element of $\basic {\textup{Conv}_m(\T)} n$. Given instead an element $B \in (\T^{N_r})^n$, the notation $M_\lambda(B)$ will denote the usual evaluation of $M_\lambda$ at $B$, seeing $M_\lambda$ as an element of $\T[x_i^{(\jj)} \mid i \in [n]; \norm \jj \le r]$, where we forget about the differential relations between the variables. We extend this notation to polynomials in $\basic {\textup{Conv}_m(\T)} n$ and $\T[x_i^{(\jj)} \mid i \in [n]; \norm \jj \le r]$, respectively. 

To make this clearer, we give a simple example here for $m=1$:
\begin{ex}
    Let $M = xx^{(3)}$, $S = 0t + 1t^3 \in \pt_{v_3}$. Taking the projection on the first four coordinates of $\Xi^{-1}_\trop(S)= (\infty, 0, \infty, 2, \infty, \dots) \in \T^\N$ we obtain the vector $B = (\infty, 0, \infty, 2) \in \T^4$. Using Definition \ref{defin:tropical-solutions}, we obtain:
    \[
        M(S) = \Phi(S) \odot \Phi(d_{v_3}^3(S)) = (1,0)\odot(0,2) = (1,2) \in \T^{(2)}.
    \]
    Evaluating $M$ at $B$ we obtain instead:
    \[
        M(B) = B_0 \odot B_3 = \infty \odot 2 = \infty \in \T.
    \]
\end{ex}

Let $V^{\trop}(E)$ denote the tropical variety of a system of tropical polynomials $E$. The following Lemma is the core of the proof of Theorem \ref{theorem:fundamental}:
\begin{lem}
\label{Lemma:truncation-is-solution}
		For every $k \in \N$, let $\pi_k \co (\T^{\N^m})^n  \rightarrow (\T^{N_k} )^n$ be the projection sending every entry to the finite vector of coordinates indexed by $\{ \jj \in \N^m \mid  \norm \jj \le J_k\}$. The following inclusion holds for every $k \in \N$: 
		\[
		 \pi_k \circ \Xi_\trop^{-1} \left (   \textup{Sol}_{\mathbf{ S}_m}(\trop_{v} (I)) \right ) \subseteq  V^{\trop} \left ( \left \{\trop_{v_K}(F_{l,\rr}) \right \}_{\substack{l=1, \dots , s \\ \norm \rr \le k}} \right ) .
		\]
	\end{lem}
        \begin{proof}
		\textbf{Step 0: Notation and set-up.} Let $S = (S_1, \dots , S_n) \in  \text{Sol}_{\mathbf S _m}(\trop_{v} (I))$ and write $S_i$ as $\sum_{\jj \in \N^m}^\infty c_{i,\jj} t^\jj$ for all $i = 1, \dots , n$.  Let 
		\[
		B_k = (B_{k,1}, \dots, B_{k,n})  :=	\pi_k \circ \Xi_\trop^{-1} (S)  \in (\T^{N_k})^n
		\]
		 where $B_{k,i} = (b_{i,\jj})_{\norm \jj \le J_k }$ with $ b_{i,\jj}:=c_{i,\jj} + v_K(\jj!)$.
        For every $\norm \rr \le k$ and $l$, let us write
        \[
                d^\rr f_l = \sum_{\lambda \in \Lambda_{l,\rr}} A_{l, \rr,\lambda} M_\lambda \in \diff{R_m} n
        \]
        for some finite set $\Lambda_{l,\rr} \subset \text{Mat}_{N_k \times n}(\N)$ and differential monomials $M_\lambda$.
	  We denote by $C_{l, \rr,\lambda}$ the valuation of $A_{l, \rr,\lambda}$, thus
	\[
                \trop_{v}(d^\rr f_l)=\sum_{\lambda \in \Lambda_{l,\rr}} C_{l, \rr,\lambda} M_\lambda \in \textup{Conv}_m(\T)[x_i^{(\jj)} \mid i \in \{1, \dots , n\}, \norm \jj \le J_k].
	\] 
        As $S \in  \text{Sol}_{\mathbf S _m}(\trop_{v} (I))$, in particular 
        \[
            S \in  \text{Sol}_{\mathbf S _m} \left (   \{ \trop_{v}(d^\rr f_l) \}_{\substack{l=1, \dots , s \\ \norm \rr \le J_k}} \right ).
        \]
        Let $\Lambda_{l,\rr}^{1} :=\{\lambda \in \Lambda_{l,\rr} \mid  \sigma_0(C_{l,\rr,\lambda}) = 1_{\textup{Conv}_m(\B)}\}$.
         The set $\Lambda_{l,\rr}^{1}$ is the support of $F_{l,\rr}$. Indeed if $\lambda \notin \Lambda_{l,\rr}^{1}$, then the $\lambda$-th term vanishes when setting $t=0$ in the definition of $F_{l,\rr}$. 
        Given that for any $\lambda \in \Lambda_{l,\rr}^1$ the element $C_{l,\rr,\lambda}$ has only one vertex (namely $0_{\N_m}$), it follows that $C_{l,\rr,\lambda} $ can be written as $ (0_{\N_m} \text{; } \alpha_{l,\rr,\lambda})$. Then, for every $\norm \rr \le k$ and $l$, we may write the tropicalization of $F_{l,\rr}$ with respect to $v_K$ as 
		\[
			\trop_{v_K}(F_{l,\rr})=\sum_{\lambda \in \Lambda_{l,\rr}^1} \alpha_{l,\rr,\lambda} M_\lambda \in \T[x_i^{(\jj)} \mid i \in \{1, \dots , n\}, \norm \jj \le J_k].
		\]
        Let now $\mathcal M :=\{\lambda \in \Lambda_{l,\rr} \mid \sigma_{0}(M_\lambda (S)) = 1_{\textup{Conv}_m(\B)}) \}$.
        
		\textbf{Step 1: If $\lambda \in \Lambda_{\lambda, \rr} \setminus \mathcal{M}$, then $M_\lambda(B_k) = \infty$.} Let $\lambda$ such that $\sigma_0(M_\lambda(S)) \neq  1_{\textup{Conv}_m(\B)}$. Then, as the product in $\textup{Conv}_m(\B)$ is given by Minkowski sum, for $\sigma_0(M_\lambda(S))$ to be strictly contained in the positive orthant it is necessary that at least one of the variables appearing in $M_\lambda$, say $x_i^{(\jj)}$, gives a weighted convex hull whose support is strictly contained in the positive orthant when evaluated at $S$. In other words, $x_i^{(\jj)}(S) \neq 1_{\textup{Conv}_m(\B)}$.
        This is as to say that the $0_{\N_m}$-th term of $d^\jj(S_i)$ is not a vertex of $\Phi_m(d^\jj(S_i))$ or, equivalently, that the $\jj$-th coefficient of $S_i$ (i.e.\ $c_{i,\jj}$) is $\infty$. Thus, the evaluation of $M_\lambda$ at $B_k$ is equal to $\infty$.

         \textbf{Step 2: For every $l = 1, \dots ,s$ and $ \norm \rr \le J_k$, $B_k$ is a solution for $\trop_{v_K}(F_{l,\rr})$.}
         From Step 1 we may assume without loss of generality that $\Lambda_{l,\rr}^{1} \cap \mathcal M \neq \varnothing$, as otherwise $B_k$ is trivially a solution for $\trop_{v_K}(F_{l,\rr})$. This implies that 
		\[
            \trop_{v}(d^\rr f_l)(S) = \bigoplus_{\lambda \in \Lambda_{l,\rr}^{1} \cap \mathcal M} (0_{\N^m} \text{; } \alpha_{l,\rr,\lambda}) \odot M_\lambda (S) = (0_{\N^m} \text{; } \alpha)
		\]
        for some $\alpha \in \R$. For every $\lambda \in \Lambda_{l, \rr}^{1} \cap \mathcal M$ we have
        \begin{align*}
				M_{\lambda}(S) &=  \bigodot_{\substack{i=1, \dots , n \\ \jj \le N_k}}  \left ( \Phi_m(d^\jj(S_i)) \right ) ^{\odot \lambda_{i,\jj}} 
                = 
                \bigodot_{\substack{i=1, \dots , n \\ \jj \le N_k}} \left ( 0_{\N^m} \text{; }  c_{i,\jj} \odot v_K(\jj!)  \right )^{\odot \lambda_{i,\jj}} = \\
				&
                = 
                \bigodot_{\substack{i=1, \dots , n \\ \jj \le N_k}} \left ( 0_{\N^m} \text{; }b_{i,\jj}  \right )^{\odot \lambda_{i,\jj}} 
                =
                \left ( 0_{\N^m} \text{; } \bigodot_{\substack{i=1, \dots , n \\ \jj \le N_k}} b_{i,\jj}^{ \odot \lambda_{i,\jj}} \right ) 
                = \left ( 0 \text{; } 
                M_{\lambda}(B_k)
                \right ) .
		\end{align*}
        As $S$ is a solution for $\trop_{v}(d^\rr f_l)$, the value $(0, \alpha)$ is attained at least twice, i.e.\ there exist two matrices of exponents $\mu, \nu \in \Lambda_{l, \rr}^{1} \cap \mathcal M$ such that 
        \[
            \left ( 0_{\N^m} \text{; } 
                \alpha_{l,r,\mu} \odot M_{\mu}(B_k)
                \right ) 
                = (0_{\N^m},\alpha) = 
                \left ( 0_{\N^m} \text{; } 
                \alpha_{l,r,\nu} \odot M_{\nu}(B_k)
                \right ) .
        \]
        This implies that for every $k \in \N$, $l = 1, \dots , s$ and $ \norm \rr \le k$, the element $B_k$ is a solution for $\trop_{v_K}(F_{l,\rr})$, proving the desired statement.
        \end{proof}

	 We are now ready to prove the fundamental theorem:

	\begin{proof}[Proof of Theorem \ref{theorem:fundamental}]
        Let $(L,v_L)$ be the extension of $(K,v_K)$, with differentially enhanced valuation, given by Corollary \ref{corolla:remark+proposition}. Fix a differential ideal $I \subset R_m^K\{x_1,\dots,x_n\}$. Let us denote by $J \subset R^L_m\{x_1, \dots x_n\}$ the differential ideal $IR^L_m\{x_1, \dots x_n\}$.
         
         The proof of the inclusion $	\trop_{\tilde v}(\text{Sol}_{R^L_m}(J)) \subseteq \text{Sol}_{\mathbf {\Gamma}^L_m}(\trop_{v_L} (J))$ is analogous to \cite[Proposition 5.2.2]{gianmereta}, thus we are left to prove the other inclusion.

        By the Fundamental Theorem of Tropical Geometry \cite[Theorem 3.2.5]{macsturm}, the following equality holds for every $k \in \N$: 
		\[
			\trop_{v_L}(A_k) = V^{\trop} \left ( \left \{\trop_{v_L}(F_{l,\rr}) \right \}_{\substack{l=1, \dots , s \\ \norm \rr \le k}} \right ) \cap (\Gamma^{N_k})^n.
		\]
		Let $S \in 	\text{Sol}_{\mathbf {\Gamma}_m^L}(\trop_{v_L} (J))$. Lemma \ref{Lemma:truncation-is-solution} applied to $L$ gives, for all $k \in \N$:
		 \[
		 \pi_k \circ \Xi_\trop^{-1}(S) \in V^{\trop} \left ( \left \{\trop_{v_L}(F_{l,\rr}) \right \}_{\substack{l=1, \dots , n \\ \norm \rr \le k}} \right ) \cap (\Gamma^{N_k})^n = \trop_{v_L}(A_k) 
		 \]
		which implies that for every $k \in \N$ there exists an element $x \in A_k$ whose tropicalization is $ \pi_k \circ \Xi_\trop^{-1}(S)$. This is equivalent to saying that $x \in (A_k)^{v_L}_S$, thus $(A_k)^{v_L}_S \neq \varnothing$ for all $k\in \N$. Since $L$ is the extension of $K$ obtained by Corollary \ref{corolla:remark+proposition} it follows that $S \in \trop_{\tilde v}(\textup {Sol}_{R^{L}_m}(J))$ 
	\end{proof}

\section{Initial forms with respect to tropical power series weight vectors} \label{section:grobner}
In this section we generalise notions introduced in \cite{FT20} (and, slightly differently, in \cite{hugao}) to the nontrivially valued case and from $m=1$ to any positive integer $m$, by defining the initial form of an element of $\diff {K(\!(t_1, \dots , t_m)\!)}{n}$ with respect to a vector of tropical power series $S \in (\T_{m,v_K})^n$. This theory is the differential counterpart to the theory of initials over a nontrivially valued field as explained for example in \cite{macsturm}. The main result of the section is the proof of the second equality of Theorem \ref{theorem:fundamental-introduction}.

Let $v_K \co K \rightarrow \T$ be a valued field: in this section we will use the same notation as in Section \ref{section:preliminaries}: $\Gamma_K$ for its value group, $R_K$ for its ring of integers, $\mathfrak m _K$ for the maximal ideal of this ring, and $k$ for its residue field. 

We recall here the definition of (an embodiment of) the field of multivariate Laurent series with integer exponents. For a more general treatement see \cite{arocacanojung, arocarond}. Given a cone $\tau$, the localization of  $K[\![\tau]\!]$ 
in the variables $t_1 , \dots , t_m$ can be realised in the following way:
\[
    K[\![\tau]\!]_{t_1 , \dots , t_m} = \bigcup_{\gamma \in \Z^m} K[\![\gamma + \tau]\!] =
    \left \{ A= \sum_{\jj \in \Z^m} a_\jj t^\jj \mid \exists \gamma \in \Z^m \text{ s.t. } \supp(A) \subseteq \gamma+\tau  \right \}
\]
Finally, choosing a vector $\omega = (\omega_1 , \dots, \omega_m) \in \R^m$ with rationally independent coordinates and denoting by $H_\omega$ the halfspace $\{v \in \R^m \mid \langle v, \omega \rangle \ge 0 \}$, define:
\[
    K(\!(t_1, \dots , t_m)\!)_\omega = \bigcup_{ \tau \subseteq H_\omega }  K[\![\tau]\!]_{t_1 , \dots , t_m}
\]
where the union is taken over all cones $\tau \subseteq H_\omega$. It is a field, and it contains $R_m$ if all the coordinates of $\omega$ are positive. From now on we will fix such a vector $\omega$ with positive, rationally independent coordinates.

\begin{oss}
Here we work with integral exponents, but it is possible (as done for instance in \cite{arocacanojung}) to consider fractional exponents as well. In that case the field obtained from the same construction as $ K(\!(t_1, \dots , t_m)\!)_\omega$ above would be an algebraically closed field containing $ K(\!(t_1, \dots , t_m)\!)_\omega$ and one could prove analogous statements as in the following for this field instead.
\end{oss}
Following similar steps as the ones in its construction, 
we will equip $K(\!(t_1, \dots , t_m)\!)_\omega$ with a valuation which is the colimit over a poset category (of extensions) of the valuations $v_\tau \co \pkcone \tau \rightarrow \textup{Conv}(\tau, \T)$. 

Let $e_i$ be the $i$-th standard basis vector for $\R^m$, for every $i \in [m]$. Given a cone $\tau$ we will denote as $\textup{Conv}(\tau, \T)_{\{e_i\} }$ the localization of the semiring $\textup{Conv}(\tau, \T)$ in the multiplicative system generated by $\{e_1 + \tau, \dots , e_m + \tau \}$. This localization can be described as the semiring whose underlying set is that of convex hulls with respect to $\tau$ of finite subsets of $\Z^m$, with vertices weighted by elements of $\T$. The operations $\oplus$ and $\odot$ in this semiring are again given by the convex hull of the union and by Minkowski sum. The semiring $\textup{Conv}(\tau, \T)_{\{e_i\}}$ is a $\T$-algebra via the same injective morphism as $\textup{Conv}(\tau, \T)$. Furthermore, for every cone $\tau \subset H_\omega$, we have a valuation \[
    \pkcone \tau _{t_1, \dots , t_m}\rightarrow \textup{Conv}(\tau, \T)_{\{e_i\}}
\]
whose definition is analogous to that of $v_\tau \co \pkcone \tau \rightarrow \textup{Conv}(\tau, \T)$ and that extends $v_\tau$ to the localisation at the variables.

Let $\mathbf{Cones}_\omega$ be the poset category of cones contained in the halfspace $H_\omega$. Sending a cone $\tau$ to the $\T$-algebra $\textup{Conv}(\tau, \T)_{\{e_i\}}$ defines a functor
\[
    \textup{Conv}(-, \T)_{\{e_i\}} \co \mathbf{Cones}_\omega \rightarrow \alg{\T}
\]
that to each inclusion $\tau \subset \rho$ associates the $\T$-algebra morphism 
\[
     \textup{Conv}(\subset, \T)_{\{e_i\}} \co \textup{Conv}(\tau, \T)_{\{e_i\}}\rightarrow  \textup{Conv}(\rho, \T)_{\{e_i\}}
\]
given by sending an element $C$ to $\sigma_0(C) + \rho$, keeping the same weights at the surviving vertices. 
\begin{ex}
We illustrate here how the morphism $ \textup{Conv}(\subset, \T)_{\{e_i\}} $ works in an explicit instance. Let $m=2$, and $\omega = (\sqrt{2},1)$. Let $\tau \subset \rho \subset H_\omega$ be respectively $\textup{pos}\langle(2,-1), (-1,2) \rangle$ and $\textup{pos}\langle (-1,2), (1,-1) \rangle$. 
Let $C = ( ((0,3);1), ((-3,5);0), ((-5,7);2) )$, depicted below.
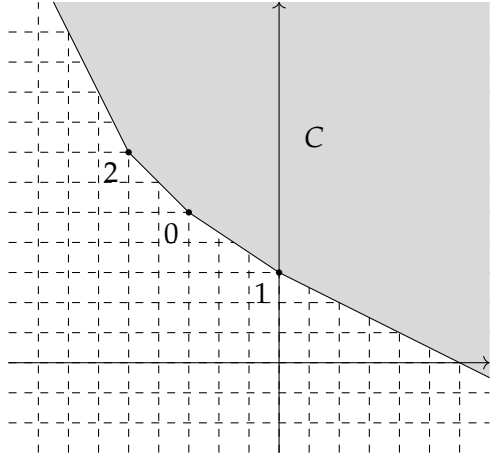
\begin{center}
\begin{tikzpicture}[scale=0.33]
\draw (-5.99,-4.99)[dashed, opacity=0.5] grid (9.99,9.99);
\fill[fill=gray!30, opacity=0.75] (-4.5,10)--(-2,5)--(0,3)--(3,1)--(10,-2.5)--(10,10);
\draw[->] (-6,-2) -- (10,-2);
\draw [->] (3,-5) -- (3,10);
\draw (-4.5,10) -- (-2,5);
\draw (-2,5)--(0,3);
\draw (0,3)--(3,1);
\draw (3,1)--(10,-2.5);
\draw (3.5,5.5) node[right] {$C$};
\draw (-2,5) node[below left] {2};
\draw (0,3) node[below left] {0};
\draw (3,1) node[below left] {1};
\fill (3,1) circle[radius=3pt];
\fill (0,3) circle[radius=3pt];
\fill (-2,5) circle[radius=3pt];
\end{tikzpicture}
\end{center}
The element $ \textup{Conv}(\subset, \T)_{\{e_i\}} (C)$ is equal to $((-5,7);2) \in \textup{Conv}(\rho, \T)_{\{e_i + \tau \mid i \in [m]\}}$:
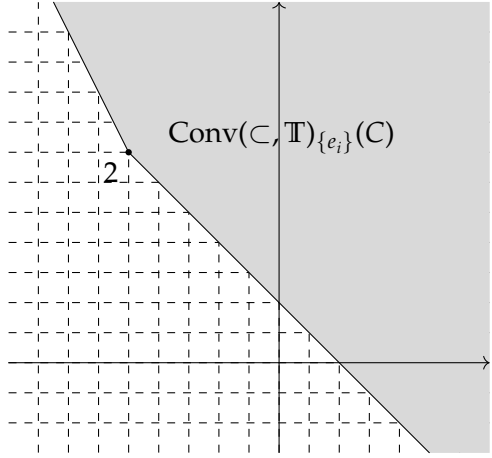
\begin{center}
\begin{tikzpicture}[scale=0.33]
\draw (-5.99,-4.99)[dashed, opacity=0.5] grid (9.99,9.99);

\fill[fill=gray!30, opacity=0.75] (-4.5,10)--(-2,5)--(8,-5)--(10,-5)--(10,10);
\draw[->] (-6,-2) -- (10,-2);
\draw [->] (3,-5) -- (3,10);
\draw (-4.5,10)--(-2,5);
\draw (-2,5)--(8,-5);
\draw (-1,5.5) node[right] {$ \textup{Conv}(\subset, \T)_{\{e_i\}} (C)$};
\draw (-2,5) node[below left] {2};
\fill (-2,5) circle[radius=3pt];
\end{tikzpicture}
\end{center}
\end{ex}
\begin{oss}
Note that, in the example above, as the image of $C$ in $\textup{Conv}(\rho, \T)_{\{e_i\}}$ is a translate of $\rho$, it is invertible. In fact, this is always true, eventually: for any cone $\tau$ and $C \in \textup{Conv}(\tau, \T)_{\{e_i\}}$, there exits a cone $\tau \subset \rho \subset H_\omega $ such that the image of $C$ in $\textup{Conv}(\rho, \T)_{\{e_i\}}$ is invertible. We will make this into a formal statement in Proposition \ref{prop:colimit-cones}.
Furthermore, by changing the vector $\omega$, we can choose any of the vertices of $C$ to be the only surviving one. Indeed, in the example above let $\omega=(1, \sqrt{2})$ then $(0,3)$ is the only vertex of $ \textup{Conv}(\subset, \T)_{\{e_i\}} (C)$ with respect to the inclusion $\tau \subset \textup{pos} \langle (3,-1), (2,-1) \rangle \subset H_\omega$. Analogously, choosing $\omega = (\sqrt{2}/2, 1)$, $(-3,5)$ is the only vertex of $ \textup{Conv}(\subset, \T)_{\{e_i\}} (C)$ with respect to the inclusion $\tau \subset \textup{pos} \langle (-1,1), (3,-2) \rangle \subset H_\omega$.
\end{oss}
\begin{prop} \label{prop:colimit-cones}
The following isomorphism of $\T$-algebras holds:
\[
    \colim_{\mathbf{Cones}_\omega} \textup{Conv}(-, \T)_{\{e_i\}}  \cong \textup{Lead}_\omega(\Z^m, \T).
\]
\end{prop}
\begin{proof}
    Firstly, the $\T$-algebra $\textup{Lead}_\omega(\Z^m, \T)$ is a co-cone for the functor $ \textup{Conv}(-, \T)_{\{e_i\}} $ via the family of $\T$-algebra morphisms $\{\mu_\tau \co \textup{Conv}(\tau, \T)_{\{e_i\}} \rightarrow \textup{Lead}_\omega(\Z^m, \T)\}_{\tau}$ defined by
    \begin{equation}\label{eq:family-cocone}
        C :=((\jj_i; \alpha_i))_{i \in [\ell]} \mapsto (\jj_{\overline i}:= \min_{\preceq_{\omega}}\{\jj_i \mid i \in [\ell]\}, \alpha_{\overline i}).
    \end{equation}
    In fact, given any two cones $\tau \subset \rho$, the diagram 
    \[
        \begin{tikzcd}
            \textup{Conv}(\tau, \T)_{\{e_i\}} 
            \arrow[rr, " \textup{Conv}(\subset{,} \T)_{\{e_i\}} "]  
            \arrow[dr, swap, "\mu_\tau"]
            & &
            \textup{Conv}(\rho, \T)_{\{e_i\}}  
            \arrow[dl, "\mu_\rho"]
            \\
             & \textup{Lead}_\omega(\Z^m, \T) & 
        \end{tikzcd}
    \]
    commutes. This is equivalent to say that given an element $C$ in $\textup{Conv}(\tau, \T)_{\{e_i\}}$, the vertex $\jj$ of $C$ attaining the minimum with respect to $\preceq_\omega$ is a vertex of the image of $C$ via any map $\textup{Conv}(\subset, \T)_{\{e_i\}} $. This indeed holds: since $\jj$ is the only vertex of $C$ lying on the hyperplane $\jj + \{x \in \R^m \mid \langle x, \omega \rangle = 0\}$ and all the others lie in $ \jj + \{x \in \R^m \mid \langle x, \omega \rangle > 0\}$, we have that $\jj$ is a vertex of the convex hull of $C$ with respect to any cone $\tau \subset \rho \subset H_\omega$.

    Thus, there is a unique morphism of $\T$-algebras
    \[
    \mu \co \colim_{\mathbf{Cones}_\omega} \textup{Conv}(-, \T)_{\{e_i\}}  \rightarrow  \textup{Lead}_\omega(\Z^m, \T).
    \]
    Since, for instance, the map $\mu_{(\R_{\ge 0})^m}$ is surjective, as the element $(\jj, \alpha)$ is the image of $C = ((\jj, \alpha)) \in  \textup{Conv}((\R_{\ge 0})^m, \T)_{\{e_i\}}$, we have that $\mu$ is surjective. For the injectivity, given two cones $\tau_1, \tau_2 \subset H_\omega$ and elements $C_1 \in \textup{Conv}(\tau_1, \T)_{\{e_i \}}$ and $C_2 \in \textup{Conv}(\tau_2, \T)_{\{e_i \}}$ such that 
    \[
        \mu_{\tau_1}(C_1) = \mu_{\tau_2}(C_2) = (\jj, \alpha) \in \textup{Lead}_\omega(\Z^m, \T), 
    \]
    by definition of the maps $\mu_{\tau_i}$ selecting the minimum of the vertices with respect to $\preceq_\omega$, there exists a cone $\tau_1, \tau_2 \subset \rho \subset H_\omega$ such that
    \[
        \textup{Conv}(\subset, \T)_{\{e_i \}}(C_1) = \textup{Conv}(\subset, \T)_{\{e_i \}}(C_2) = ((\jj, \alpha)) \in \textup{Conv}(\rho, \T)_{\{e_i \}}
    \]
    thus $C_1$ and $C_2$ are identified in the colimit, and the map $\mu$ is injective, thus giving the desired isomorphism.
\end{proof}
The semiring $\textup{Lead}_\omega(\Z^m, \T)$ is a semifield and by Proposition \ref{prop:colimit-cones} above we have a valuation 
\begin{equation} \label{eq:v-omega}
    v_\omega \co K(\!(t_1, \dots , t_m)\!)_\omega \rightarrow \textup{Lead}_\omega(\Z^m, \T).
\end{equation}
defined by sending a series $A=\sum_{\jj \in \Z^m} a_\jj t^\jj$ to $(\overline \jj:=\min_{\preceq_\omega}\supp( A), v_K(a_{\overline \jj} ))$, that becomes surjective restricting the weights to the value group $\Gamma_K$ of $v_K$. Notice that the restriction of $v_\omega$ to $\pkcone \tau$ for any cone $\tau \subset H_\omega$, is the valuation $v_{\tau, \omega}$ of Example \ref{ex:omega-val-tau}. 

If the valuation $v_K \co K \rightarrow \T$ has a section then the corresponding extended valuation $v_\omega \co K(\!(t_1, \dots , t_m)\!)_\omega \rightarrow \textup{Lead}_\omega(\Z^m, \widetilde{\Gamma}_K)$ has a section as well. Explicitly, let $\pi_K$ be a uniformizer for $v_K$, then the map 
\begin{align*}
    \varphi_\omega \co \textup{Lead}_\omega(\Z^m, \widetilde{\Gamma}_K) &\rightarrow K(\!(t_1, \dots t_m)\!)_\omega \\
   (\jj,\alpha) &\mapsto \pi_K^{\alpha} t^{\jj} 
\end{align*}
is a section for the valuation $v_\omega$, and we will drop the subscript $\omega$ when it is clear from the context. Given a cone $\tau \subset H_\omega$, let again $v_\tau$ be the valuation 
\[
    v_\tau \co \pkcone \tau _{t_1, \dots t_m} \rightarrow \textup{Conv}(\tau, \T)_{\{e_i \}}
\]
extending the valuation $v_\tau$ of Example \ref{ex:conv-vals}. It is easy to check that the valuation ring with respect to $v_\tau$ and the maximal ideal of this ring are the following: 
\begin{align*}
    R_{v_\tau} :&= \{A \in \pkcone \tau _{t_1, \dots t_m} \mid v_\tau(A) \preceq (0_{\Z^m},0)\} = 
    \left \{A =\sum_{\jj \in \Z^m} a_\jj t^\jj \in \pkcone \tau \mid a_{0_{\Z^m}} \in R_K \right  \}, 
    \\
    \mathfrak m _{v_\tau} :&= \{A \in \pkcone \tau _{t_1, \dots t_m} \mid v_\tau(A) \prec (0_{\Z^m},0)\} =
    \left \{A =\sum_{\jj \in \Z^m} a_\jj t^\jj \in \pkcone \tau \mid  a_{0_{\Z^m}} \in \mathfrak m_K \right \}, 
\end{align*}   
and the residue field $ R_{v_\tau}/\mathfrak m _{v_\tau}$ is isomorphic to $k=R_K/\mathfrak m_K$.
From the equalities above, and from the definition of the field $K(\!(t_1, \dots t_m)\!)_\omega$, the following hold:
\begin{align*}
    R_{v_\omega} :&= \{A \in K(\!(t_1, \dots t_m)\!)_\omega \mid v_\omega(A) \preceq (0_{\Z^m},0)\} = \bigcup_{\tau \subset H_\omega} R_{v_\tau} =
    \\
    &= \left \{A \in K(\!(t_1, \dots t_m)\!)_\omega \mid  A \in R_{v_\tau} \text{ for some cone } \tau \subset H_\omega \right \} = 
    \\
    &= \left \{ A=\sum_{\jj \in \Z^m} a_\jj t^\jj \in K(\!(t_1, \dots t_m)\!)_\omega \mid  \supp(A) \subset H_\omega \text{ and } a_{0_{\Z^m}} \in R_K \right \} ;
    \\
    \mathfrak m _{v_\omega} :&= \{A \in K(\!(t_1, \dots , t_m)\!)_\omega \mid v_\omega(A) \prec (0_{\Z^m},0)\} = \bigcup_{\tau \subset H_\omega} \mathfrak m_{v_\tau} =
    \\
    &= \left \{A \in K(\!(t_1, \dots t_m)\!)_\omega \mid A \in \mathfrak m_{v_\tau} \text{ for some cone } \tau \subset H_\omega \right \} =
    \\
    &= \left \{A=\sum_{\jj \in \Z^m} a_\jj t^\jj \in K(\!(t_1, \dots t_m)\!)_\omega \mid  \supp(A) \subset H_\omega \text{ and } a_{0_{\Z^m}} \in \mathfrak m_K \right \}.
\end{align*}   
The residue field $R_{v_\omega}/\mathfrak m _{v_\omega}$ is isomorphic to $k$. We will write $\overline A \in k$ for the residue class of an element $A \in R_{v_\omega}$ and if $f \in \diff {K(\!(t_1 , \dots , t_m)\!)_\omega}{n}$ has coefficients in $R_{v_\omega}$, let $\overline f \in \diff k n$ be the differential polynomial obtained by taking the residue class of its coefficients.

Given $S := (S_1, \dots , S_n) \in (\T_{m,v_K})^n$, the elements $\Phi_m(d^\jj S_i)$ and $\trop_v(f)(S)$ live in $\textup{Conv}_m(\T) \subset \textup{Conv}_m(\T)_{\{e_i\}}$, thus we can look at their image via the map $\mu_m :=\mu_{(\R_{\ge 0})^m}$ of \ref{eq:family-cocone}. For every $\jj \in \N^m$ such that $\Phi_m(d^\jj S_i)\neq \varnothing$ let 
\[
(\fatg_{i, \jj}, \alpha_{i, \jj}) := \mu_m(\Phi_m(d^\jj S_i)).
\]
In the same way, if $\trop_v(f)(S) \neq \varnothing$, let $(\hh, \alpha) :=\mu_m( \trop_v(f)(S))$. 

\begin{defin} \label{definition:initial}
    Let $S \in (\T_{m,v_K})^n$ and $f \in \diff {K(\!(t_1, \dots , t_m)\!)_\omega}{n}$. Writing $f = \sum_{\lambda \in \Lambda} A_\lambda \prod_{i,\jj}(x_i^{(\jj)})^{\lambda_{i,\jj}}$ as in Section \ref{section:theorem}, we define:
    \[
    h_{(S,\omega)} := \varphi  \left ( (-\hh; -\alpha) \right ) \sum_{\lambda \in \Lambda} A_\lambda \prod_{i,\jj} \left ( \varphi  \left ( ( \fatg_{i,\jj}, \alpha_{i,\jj} ) \right )  x_i^{(\jj)} \right ) ^{\lambda_{i,\jj}} \in \diff {K(\!(t_1, \dots , t_m)\!)_\omega} n.
    \] 
    The differential polynomial $h_{(S,\omega)}$ has coefficients in $R_{v_\omega}$.
    The initial form $\initial_{(S,\omega)}(f) \in \diff k n$ of $f$ with respect to $(S, \omega)$ is defined as follows:
    \[
    \initial_{(S, \omega)}(f) :=
        \begin{cases}
            \overline {h_{(S, \omega)}}   & \text{if $\trop_v(f)(S) \neq \infty$} \\
            & \\
            0 & \text{if $\trop_v(f)(S) = \infty$}
        \end{cases}
    \]
\end{defin}

\begin{oss}\label{oss:initial}
    There is at least one coefficient in the differential polynomial $h_{(S,\omega)}$ above with valuation $(0_{\Z^m},0)=1_{\textup{Lead}_\omega(\Z^m,\T)}$ with respect to $v_\omega$, thus $\initial_{(S, \omega)}(f) = 0$ if and only if $\trop_v(f)(S) = \varnothing \in \textup{Conv}_m(\T)$. Indeed
    \begin{align*}
            h_{(S, \omega)}  :&= \varphi  \left ( (-\hh; -\alpha) \right ) \sum_{\lambda \in \Lambda} A_\lambda \prod_{i,\jj} \left ( \varphi  \left ( ( \fatg_{i,\jj}, \alpha_{i,\jj} ) \right )  x_i^{(\jj)} \right ) ^{\lambda_{i,\jj}} = \\
            & = \varphi  \left ( (-\hh, - \alpha)  \right ) \sum_{\lambda \in \Lambda} A_\lambda \prod_{i,\jj}  \varphi  \left ( (\fatg_{i,\jj}, \alpha_{i,\jj})  \right )   ^{\lambda_{i,\jj}}  x^\lambda
	\end{align*}
    thus the $v_\omega$ valuation of the coefficient of $x^\lambda$ in $h_{(S,\omega)}$ is
    \[
    (-\hh, -\alpha) \odot v_\omega(A_\lambda) \odot \bigodot_{i,\jj} (\fatg_{i,\jj}, \alpha_{i,\jj})^{\odot \lambda_{i,\jj}}. 
    \]
    Given that $(\hh, \alpha) = \mu_m( \trop_v(f)(S))$ is the minimum of the elements $v_\omega(A_\lambda) \odot \bigodot_{i,\jj} (\fatg_{i,\jj}, \alpha_{i,\jj})^{\odot \lambda_{i,\jj}}$, this proves that the coefficients of the polynomial $h_{(S, \omega)}$ are elements of $R_{v_\omega}$. In particular, the only coefficients with valuation $(0_{\Z^m},0)$ are those of the monomials where $(\hh, \alpha)$ is attained.
\end{oss}
\begin{oss}
    Notice that for $m=1$, the initial form of a differential polynomial, and so the initial ideal of a differential ideal (see Definition \ref{def:initial-ideal}) depends only on $S$: in this case in fact, up to scaling, there is a unique choice for $\omega$.
\end{oss}

\begin{defin} \label{def:initial-ideal}
    Let $I \subseteq \diff {R_m^K} n$ be a differential ideal. Then the initial ideal of $I$ with respect to $(S, \omega)$ is the (algebraic) ideal: 
    \[
        \Initial_{(S, \omega)}(I):=\langle \initial_{(S, \omega)}(f) \mid f \in I \rangle \subseteq \diff k n
    \]
\end{defin}

\begin{oss}
Similarly to what already noticed in \cite[Example 3.5]{FT20}, given $f \in I$ in general we have 
\[d(\initial_{(S,\omega)}(f)) \neq \initial_{(S, \omega)}(df).\]
\end{oss}
We now generalize \cite[Lemma 2.6]{hugao} and \cite[Lemma 3.8]{FT20} to the nontrivially valued multivariate setting:
\begin{lem} \label{lemma:initial}
    Let $S \in (\T_{m,v_K})^n$ and $I \subseteq \diff {R^K_m} n$ be a differential ideal, for every $F \in \Initial_{(S, \omega)}(I)$ there exists $f \in I$ such that $ F = \initial_{(S, \omega)}(f)$. 
\end{lem}
\begin{proof}
    As in Definition \ref{definition:initial}, for all $\jj \in \N^m$ such that $\Phi_m(d^\jj S_i)\neq \varnothing$, let 
    \[
    (\fatg_{i,\jj}, \alpha_{i,\jj}) := \mu_m(\Phi_m(d^\jj S_i)).
    \]
    Given $F \in \Initial_{(S, \omega)}(I)$, we can write it as  $\sum_{\lambda \in \Lambda} a_\lambda x^\lambda \initial_{(S, \omega)}(f_\lambda)$ for some $a_\lambda \in k$ and $f_\lambda \in I$. Let 
    \[
        (\hh_\lambda, \alpha_\lambda):=\mu_{m}(\trop_v(f_\lambda)(S)) \quad \quad 
        (\hh'_\lambda, \alpha'_\lambda):= (\hh_\lambda, \alpha_\lambda) \odot \bigodot_{i,\jj}  (\fatg_{i,\jj}, \alpha_{i,\jj})^{\odot\lambda_{i,\jj}}
    \]
    and 
    \[
        f := \sum_{\lambda \in \Lambda} A_\lambda \varphi \left ((-\hh'_\lambda, -\alpha'_\lambda) \right ) x^\lambda f_\lambda
    \]
    for elements $A_\lambda \in R_m$ such that $v(A_\lambda) = v_\omega(A_\lambda) = (0_{\Z^m},0)$ and $\overline {A_\lambda} = a_\lambda$: notice that $\mu_m(\trop_{v}(f)(S)) = (0_{\Z^m},0)$.
    
    We claim that the initial form of $f$ with respect to $(S, \omega)$ is $F$. Let us compute $\initial_{(S,\omega)}(f)$:
    \begin{align*}
    \initial_{(S,\omega)}(f) & = \overline {h_{(S,\omega)}} = \\
    & = \overline { \sum_{\lambda \in \Lambda} A_\lambda \varphi \left ((-\hh'_\lambda, -\alpha'_\lambda) \right ) \prod_{i,\jj} \left ( \varphi  \left ( (\fatg_{i,\jj}, \alpha_{i,\jj})  \right )  x_i^{(\jj)} \right) ^{\lambda_{i,\jj}} f_\lambda \left ( \varphi \left ( (\fatg_{i,\jj}, \alpha_{i,\jj})  \right )  x_i^{(\jj)} \right)} = \\
    & =  \overline{ \sum_{\lambda \in \Lambda} A_\lambda \varphi \left ((-\hh'_\lambda, -\alpha'_\lambda) \right ) \prod_{i,\jj} \left ( \varphi  \left ( (\fatg_{i,\jj}, \alpha_{i,\jj})  \right )  \right ) ^{\lambda_{i,\jj}}  x^\lambda f_\lambda \left ( \varphi \left ( (\fatg_{i,\jj}, \alpha_{i,\jj})  \right )  x_i^{(\jj)} \right ) }= \\
    & = \sum_{\lambda \in \Lambda} \overline{A_\lambda  \varphi \left ((-\hh_\lambda, -\alpha_\lambda) \right )  f_\lambda \left ( \varphi \left ((\fatg_{i,\jj}, \alpha_{i,\jj})  \right )  x_i^{(\jj)} \right)  x^\lambda }   .\\
    \end{align*}
    As $\overline{A_\lambda} = a_\lambda$ by definition and $\overline{ \varphi \left ((-\hh_\lambda, -\alpha_\lambda) \right ) f_\lambda \left ( \varphi \left ((\fatg_{i,\jj}, \alpha_{i,\jj})  \right )  x_i^{(\jj)} \right) }$ is the initial form of $f_\lambda$ with respect to $(S,\omega)$, we have
    \[
        \initial_{(S, \omega)}(f) = \sum_{\lambda \in \Lambda} a_\lambda x^\lambda \initial_{(S, \omega)}(f_\lambda) = F
    \]
    that proves the claim.
\end{proof}
Since given any vector $\omega \in (\R_{>0})^m$ with rationally independent coordinates we have that $R^K_m \subset K(\!(t_1, \dots t_m)\!)_\omega$, we can consider how the initial of $I$ with respect to $(S,\omega)$ varies by varying the vector $\omega$. Let 
\[
     \Sigma_\omega :=\{S \in \widetilde{\Gamma}_K[\![t_1, \dots t_m]\!]^n_{v_K} \mid \Initial_{(S, \omega)}(I) \text{ does not contain a monomial}\} 
\]
then we have:

\begin{teorema}\label{theorem:fundamental-initial}
Let $(K,v_K)$ be an algebraically closed valued field of characteristic 0 and $\mathbf{\Gamma}_m^K$ be the tropical pair as in Remark \ref{oss:restr-of-diff-enh}. Let $I \subset R^K_m\{x_1,\dots,x_n\}$ be a differential ideal, then the following equality holds: 
\begin{equation}\label{eq:solutions-initial}
         \textup{Sol}_{\mathbf \Gamma^K_m}(\trop_{v} (I)) 
         =
        \bigcap_{\substack{\omega \in (\R_{>0})^m
                \\ 
                \text{with $\Q$-ind. coords}
                }} \Sigma_\omega
\end{equation}
\end{teorema}
\begin{oss}
    In contrast to the statement of Theorem \ref{theorem:fundamental}, here we do not need to extend the field (and thus possibly its value group) as the equality of Equation \ref{eq:solutions-initial} happens in the tropical world and the extension of $K$ was needed to deal with issues arising from the fibers of the tropicalization map.
\end{oss}
\begin{proof}[Proof of Theorem \ref{theorem:fundamental-initial}]

        Let us start from the "$\subseteq$" inclusion. Given $S \in \textup{Sol}_{\mathbf \Gamma_m}(\trop_{v} (I))$, by definition, for every $f \in I$, either $\trop_v(f)(S)=\varnothing$ or each of its vertices is a vertex of the evaluation at $S$ of at least two monomials, with the same weight it carries in $\trop_v(f)(S)$. 
        
        For any vertex $\jj$ of $\trop_v(f)(S)$, let $\alpha$ be its weight and choose $\omega_\jj \in (\R_{>0})^m$ with rationally independent coordinates, such that the image of $\trop_v(f)(S)$ via the map
        \[
            \mu_{m,\jj} \co \textup{Conv}_m(\T)_{\{e_i\}} \rightarrow \textup{Lead}_\omega(\Z^m, \T)
        \]
        is equal to $(\jj, \alpha)$, i.e.\ $\jj$ minimizes $\langle \jj, \omega \rangle$ among the vertices of $\trop_v(f)(S)$.
        From Remark \ref{oss:initial} we know that for every $f \in I$ the initial $\initial_{(S, \omega)}(f)$ is either 0 or it is at least a binomial. The ideal $\Initial_{(S, \omega)}(I)$ cannot contain any monomial: indeed if $F \in \Initial_{(S, \omega)}(I)$ is a monomial, then from Lemma \ref{lemma:initial}, there exists $f \in I$ such that $\initial_{(S, \omega)}(f) = F$ is a monomial, and that contradicts the fact that $S$ is a solution for $\trop_{v} (I)$.
        
        For what concerns the opposite inclusion: assume $S$ is not in $\textup{Sol}_{\mathbf \Gamma_m}(\trop_{v} (I))$. Then there would exist a polynomial $f \in I$ and a weighted vertex $(\jj, \alpha)$ of $\trop_v(f)(S)$ that appears as a vertex of the evaluation at $S$ of exactly one of the monomials of $f$. Thus, choosing $\omega_\jj$ such that $\jj$ minimizes $\langle \jj, \omega \rangle$ among the vertices of $\trop_v(f)(S)$ as above, thanks to the discussion of Remark \ref{oss:initial}, we have that $\initial_{(S, \omega)}(f)$ is a monomial. Thus $S$ does not belong to the right hand side of \ref{eq:solutions-initial}. This concludes the proof.
\end{proof}

\begin{oss}
The theory developed above presents interesting links with some objects and results of \cite{initialslara} and \cite{kalina}, that the second author is investigating further. 
\end{oss}

\begin{ex}
    Let $m=2$, $n=1$ and $p=2$. Let $\C_2[\![t_1, t_2]\!]$ be endowed with the differentially enhanced valuation:
    \[
    \begin{tikzcd}
        & \T[\![t_1, t_2]\!] _{v_2} \arrow[d, "\Phi_2"] \\
        \C_2[\![t_1, t_2]\!]  \arrow[ur, "\widetilde{v}"] \arrow[r,swap,"v"] & \textup{Conv}_2(\T) 
    \end{tikzcd}
    \]
as in point (4) of Example \ref{ex:diff-enhancements}. Consider the following differential polynomial in $\C_2[\![t_1, t_2]\!]\{x\}$: 
    \[
        F:=\left(6t_1+\frac{5}{2}t_2 \right) (x^{(0,1)})^2 + 3t_1t_2x^{(2,3)} + (4t_1t_2^3 + t_1^2t_2 - 4t_1^3)x^{(1,0)}.
    \]
The tropicalization of $F$ is 
\[\trop_v(F)=    
    \begin{tikzpicture}        
        [line cap=round,line    join=round,x=1cm,y=1cm,baseline={0.5cm-0.5*height("$=$")}, scale=0.32]
        \draw (-0.99,-0.99)[dashed, opacity=0.5] grid (3.99,3.99);
        \fill[fill=gray!30, opacity=0.75] (0,4)--(0,1)--(1,0)--(4,0)--(4,4);
        \draw[->] (-1,0) -- (4,0);
        \draw [->] (0,-1) -- (0,4);
        \draw (0,4) -- (0,1);
        \draw (0,1)--(1,0);
        \draw (1,0)--(4,0);
        \draw (1,0) node[below] {\footnotesize 1};
        \draw (0,1) node[left]{\footnotesize -1};
        \fill (1,0) circle[radius=3pt];
        \fill (0,1) circle[radius=3pt];
    \end{tikzpicture}
        (x^{(0,1)})^2 +
    \begin{tikzpicture}        
        [line cap=round,line    join=round,x=1cm,y=1cm,baseline={0.5cm -0.5*height("$=$")}, scale=0.32]
        \draw (-0.99,-0.99)[dashed, opacity=0.5] grid (3.99,3.99);
        \fill[fill=gray!30, opacity=0.75] (1,4)--(1,1)--(4,1)--(4,4);
        \draw[->] (-1,0) -- (4,0);
        \draw [->] (0,-1) -- (0,4);
        \draw (1,4) -- (1,1);
        \draw (1,1)--(4,1);
        \draw (1.2,1.2) node[below left ] {\footnotesize 0};
        \fill (1,1) circle[radius=3pt];
    \end{tikzpicture}
     x^{(2,3)} +
    \begin{tikzpicture}        
         [line cap=round,line    join=round,x=1cm,y=1cm,baseline={0.5cm-0.5*height("$=$")}, scale=0.32]
        \draw (-0.99,-0.99)[dashed, opacity=0.5] grid (3.99,3.99);
        \fill[fill=gray!30, opacity=0.75] (1,4)--(1,3)--(2,1)--(3,0)--(4,0)--(4,4);
        \draw[->] (-1,0) -- (4,0);
        \draw [->] (0,-1) -- (0,4);
        \draw (1,4) -- (1,3);
        \draw (1,3)--(2,1);
        \draw (2,1)--(3,0);
        \draw (3,0)--(4,0);
        \draw (1,3) node[below left] { \footnotesize 2};
        \draw (2.2,1.2) node[below left]{ \footnotesize 0};
        \draw (3,0) node[below]{ \footnotesize 2};
        \fill (3,0) circle[radius=3pt];
        \fill (1,3) circle[radius=3pt];
        \fill (2,1) circle[radius=3pt];
    \end{tikzpicture}
    x^{(1,0)} \in \textup{Conv}_2(\T)\{x\}_{\textit{basic}}.
\]
We claim that the tropical power series
\[ 
    S = 1t_1 + 1t_1t_2 +  1t_1^2t_2^5 + (-1)t_1^3t_2^3 \in  \T[\![t_1, t_2]\!] _{v_2}
\]
is a solution for the tropicalization of $F$. Indeed we have: 
\begin{align*}\trop_v(F)(S) & =    
    \begin{tikzpicture}        
        [line cap=round,line  join=round,x=1cm,y=1cm,baseline={0.5cm-0.5*height("$=$")}, scale=0.32]
        \draw (-0.99,-0.99)[dashed, opacity=0.5] grid (3.99,3.99);
        \fill[fill=gray!30, opacity=0.75] (0,4)--(0,1)--(1,0)--(4,0)--(4,4);
        \draw[->] (-1,0) -- (4,0);
        \draw [->] (0,-1) -- (0,4);
        \draw (0,4) -- (0,1);
        \draw (0,1)--(1,0);
        \draw (1,0)--(4,0);
        \draw (1,0) node[below] {\footnotesize 1};
        \draw (0,1) node[left]{\footnotesize -1};
        \fill (1,0) circle[radius=3pt];
        \fill (0,1) circle[radius=3pt];
    \end{tikzpicture}
    \odot (\Phi_2(d^{(0,1)}S))^2 \oplus
    \begin{tikzpicture}        
        [line cap=round,line    join=round,x=1cm,y=1cm,baseline={0.5cm -0.5*height("$=$")}, scale=0.32]
        \draw (-0.99,-0.99)[dashed, opacity=0.5] grid (3.99,3.99);
        \fill[fill=gray!30, opacity=0.75] (1,4)--(1,1)--(4,1)--(4,4);
        \draw[->] (-1,0) -- (4,0);
        \draw [->] (0,-1) -- (0,4);
        \draw (1,4) -- (1,1);
        \draw (1,1)--(4,1);
        \draw (1.2,1.2) node[below left ] {\footnotesize 0};
        \fill (1,1) circle[radius=3pt];
    \end{tikzpicture}
      \odot \Phi_2(d^{(2,3)}S) \oplus
    \begin{tikzpicture}        
         [line cap=round,line    join=round,x=1cm,y=1cm,baseline={0.5cm-0.5*height("$=$")}, scale=0.32]
        \draw (-0.99,-0.99)[dashed, opacity=0.5] grid (3.99,3.99);
        \fill[fill=gray!30, opacity=0.75] (1,4)--(1,3)--(2,1)--(3,0)--(4,0)--(4,4);
        \draw[->] (-1,0) -- (4,0);
        \draw [->] (0,-1) -- (0,4);
        \draw (1,4) -- (1,3);
        \draw (1,3)--(2,1);
        \draw (2,1)--(3,0);
        \draw (3,0)--(4,0);
        \draw (1,3) node[below left] { \footnotesize 2};
        \draw (2.2,1.2) node[below left]{ \footnotesize 0};
        \draw (3,0) node[below]{ \footnotesize 2};
        \fill (3,0) circle[radius=3pt];
        \fill (1,3) circle[radius=3pt];
        \fill (2,1) circle[radius=3pt];
    \end{tikzpicture}
    \odot \Phi_2(d^{(1,0)}S)  = \\ 
& =    
    \begin{tikzpicture}        
        [line cap=round,line    join=round,x=1cm,y=1cm,baseline={0.5cm-0.5*height("$=$")}, scale=0.32]
        \draw (-0.99,-0.99)[dashed, opacity=0.5] grid (3.99,3.99);
        \fill[fill=gray!30, opacity=0.75] (0,4)--(0,1)--(1,0)--(4,0)--(4,4);
        \draw[->] (-1,0) -- (4,0);
        \draw [->] (0,-1) -- (0,4);
        \draw (0,4) -- (0,1);
        \draw (0,1)--(1,0);
        \draw (1,0)--(4,0);
        \draw (1,0) node[below] {\footnotesize 1};
        \draw (0,1) node[left]{\footnotesize -1};
        \fill (1,0) circle[radius=3pt];
        \fill (0,1) circle[radius=3pt];
    \end{tikzpicture}
    \odot 
    \begin{tikzpicture}        
        [line cap=round,line    join=round,x=1cm,y=1cm,baseline={0.5cm-0.5*height("$=$")}, scale=0.32]
        \draw (-0.99,-0.99)[dashed, opacity=0.5] grid (3.99,3.99);
        \fill[fill=gray!30, opacity=0.75] (2,4)--(2,0)--(4,0)--(4,4);
        \draw[->] (-1,0) -- (4,0);
        \draw [->] (0,-1) -- (0,4);
        \draw (2,4) -- (2,0);
        \draw (2,0)--(2,4);
        \draw (2,0) node[below] {\footnotesize 2};
        \fill (2,0) circle[radius=3pt];
    \end{tikzpicture}
    \oplus
    \begin{tikzpicture}        
        [line cap=round,line    join=round,x=1cm,y=1cm,baseline={0.5cm -0.5*height("$=$")}, scale=0.32]
        \draw (-0.99,-0.99)[dashed, opacity=0.5] grid (3.99,3.99);
        \fill[fill=gray!30, opacity=0.75] (1,4)--(1,1)--(4,1)--(4,4);
        \draw[->] (-1,0) -- (4,0);
        \draw [->] (0,-1) -- (0,4);
        \draw (1,4) -- (1,1);
        \draw (1,1)--(4,1);
        \draw (1.2,1.2) node[below left ] {\footnotesize 0};
        \fill (1,1) circle[radius=3pt];
    \end{tikzpicture}
      \odot 
      \begin{tikzpicture}        
        [line cap=round,line    join=round,x=1cm,y=1cm,baseline={0.5cm-0.5*height("$=$")}, scale=0.32]
        \draw (-0.99,-0.99)[dashed, opacity=0.5] grid (3.99,3.99);
        \fill[fill=gray!30, opacity=0.75] (0,4)--(0,2)--(1,0)--(4,0)--(4,4);
        \draw[->] (-1,0) -- (4,0);
        \draw [->] (0,-1) -- (0,4);
        \draw (0,4) -- (0,2);
        \draw (0,2)--(1,0);
        \draw (1,0)--(4,0);
        \draw (1,0) node[below] {\footnotesize 1};
        \draw (0,2) node[left]{\footnotesize 3};
        \fill (1,0) circle[radius=3pt];
        \fill (0,2) circle[radius=3pt];
    \end{tikzpicture}
    \oplus
    \begin{tikzpicture}        
         [line cap=round,line    join=round,x=1cm,y=1cm,baseline={0.5cm-0.5*height("$=$")}, scale=0.32]
        \draw (-0.99,-0.99)[dashed, opacity=0.5] grid (3.99,3.99);
        \fill[fill=gray!30, opacity=0.75] (1,4)--(1,3)--(2,1)--(3,0)--(4,0)--(4,4);
        \draw[->] (-1,0) -- (4,0);
        \draw [->] (0,-1) -- (0,4);
        \draw (1,4) -- (1,3);
        \draw (1,3)--(2,1);
        \draw (2,1)--(3,0);
        \draw (3,0)--(4,0);
        \draw (1,3) node[below left] { \footnotesize 2};
        \draw (2.2,1.2) node[below left]{ \footnotesize 0};
        \draw (3,0) node[below]{ \footnotesize 2};
        \fill (3,0) circle[radius=3pt];
        \fill (1,3) circle[radius=3pt];
        \fill (2,1) circle[radius=3pt];
    \end{tikzpicture}
    \odot 
    \begin{tikzpicture}        
        [line cap=round,line    join=round,x=1cm,y=1cm,baseline={0.5cm-0.5*height("$=$")}, scale=0.32]
        \draw (-0.99,-0.99)[dashed, opacity=0.5] grid (3.99,3.99);
        \fill[fill=gray!30, opacity=0.75] (0,4)--(0,0)--(4,0)--(4,4);
        \draw[->] (-1,0) -- (4,0);
        \draw [->] (0,-1) -- (0,4);
        \draw (0,0) node[below left ] {\footnotesize 1};
        \fill (0,0) circle[radius=3pt];
    \end{tikzpicture}
    = \\
    & =    
    \begin{tikzpicture}        
        [line cap=round,line    join=round,x=1cm,y=1cm,baseline={0.5cm-0.5*height("$=$")}, scale=0.32]
        \draw (-0.99,-0.99)[dashed, opacity=0.5] grid (3.99,3.99);
        \fill[fill=gray!30, opacity=0.75] (2,4)--(2,1)--(3,0)--(4,0)--(4,4);
        \draw[->] (-1,0) -- (4,0);
        \draw [->] (0,-1) -- (0,4);
        \draw (2,4) -- (2,1);
        \draw (2,1)--(3,0);
        \draw (3,0)--(4,0);
        \draw (2,1) node[left] {\footnotesize 1};
        \draw (3,0) node[below]{\footnotesize 3};
        \fill (3,0) circle[radius=3pt];
        \fill (2,1) circle[radius=3pt];
    \end{tikzpicture}
    \oplus
    \begin{tikzpicture}        
        [line cap=round,line    join=round,x=1cm,y=1cm,baseline={0.5cm -0.5*height("$=$")}, scale=0.32]
        \draw (-0.99,-0.99)[dashed, opacity=0.5] grid (3.99,3.99);
        \fill[fill=gray!30, opacity=0.75] (1,4)--(1,3)--(2,1)--(4,1)--(4,4);
        \draw[->] (-1,0) -- (4,0);
        \draw [->] (0,-1) -- (0,4);
        \draw (1,4) -- (1,3);
        \draw (1,3)--(2,1);
        \draw (2,1)--(4,1);
        \draw (1,3) node[below left ] {\footnotesize 3};
        \draw (2.2,1.2) node[below left ] {\footnotesize 1};
        \fill (1,3) circle[radius=3pt];
        \fill (2,1) circle[radius=3pt];
    \end{tikzpicture}
    \oplus
    \begin{tikzpicture}        
         [line cap=round,line    join=round,x=1cm,y=1cm,baseline={0.5cm-0.5*height("$=$")}, scale=0.32]
        \draw (-0.99,-0.99)[dashed, opacity=0.5] grid (3.99,3.99);
        \fill[fill=gray!30, opacity=0.75] (1,4)--(1,3)--(2,1)--(3,0)--(4,0)--(4,4);
        \draw[->] (-1,0) -- (4,0);
        \draw [->] (0,-1) -- (0,4);
        \draw (1,4) -- (1,3);
        \draw (1,3)--(2,1);
        \draw (2,1)--(3,0);
        \draw (3,0)--(4,0);
        \draw (1,3) node[below left] { \footnotesize 3};
        \draw (2.2,1.2) node[below left]{ \footnotesize 1};
        \draw (3,0) node[below]{ \footnotesize 3};
        \fill (3,0) circle[radius=3pt];
        \fill (1,3) circle[radius=3pt];
        \fill (2,1) circle[radius=3pt];
    \end{tikzpicture}
    = \\
    & = 
    \begin{tikzpicture}        
         [line cap=round,line    join=round,x=1cm,y=1cm,baseline={0.5cm-0.5*height("$=$")}, scale=0.32]
        \draw (-0.99,-0.99)[dashed, opacity=0.5] grid (3.99,3.99);
        \fill[fill=gray!30, opacity=0.75] (1,4)--(1,3)--(2,1)--(3,0)--(4,0)--(4,4);
        \draw[->] (-1,0) -- (4,0);
        \draw [->] (0,-1) -- (0,4);
        \draw (1,4) -- (1,3);
        \draw (1,3)--(2,1);
        \draw (2,1)--(3,0);
        \draw (3,0)--(4,0);
        \draw (1,3) node[below left] { \footnotesize 3};
        \draw (2.2,1.2) node[below left]{ \footnotesize 1};
        \draw (3,0) node[below]{ \footnotesize 3};
        \fill (3,0) circle[radius=3pt];
        \fill (1,3) circle[radius=3pt];
        \fill (2,1) circle[radius=3pt];
    \end{tikzpicture}
\end{align*}
The sum above tropically vanishes, thus $S$ is a solution for $\trop_v(F)$. We now show that the initial of $F$ with respect to $(S,\omega)$ is not a monomial for every $\omega$. The vertices of $\trop_v(F)(S)$ are $V:=\{(3,0),(2,1),(1,3)\}$. Let $\omega$ with positive, rationally independent coordinates such that $(3,0) = \min_{\preceq_\omega} V$. Then we have: 
\begin{align*}
    h_{(S,\omega)} & = 2^{-3}t_1^{-3} \left ( (6t_1+\frac{5}{2}t_2 ) \cdot (2t_1x^{(0,1)})^2 + 3t_1t_2 \cdot 2t_1x^{(2,3)} + (4t_1t_2^3 + t_1^2t_2 - 4t_1^3) \cdot 2x^{(1,0)}
    \right) = \\
    & = 3+\frac{5}{4}t_2t_1^{-1}  (x^{(0,1)})^2 + \frac{3}{4}t_1^{-1}t_2 x^{(2,3)} + (t_1^{-2}t_2^3 + \frac{1}{4}t_1^{-1}t_2 - 1) x^{(1,0)}
\end{align*}
and $\initial_{(S, \omega)}(F)$ is the reduction of $ h_{(S,\omega)}$ modulo $\mathfrak m_{v_\omega}$:
\[
    \initial_{(S, \omega)}(F)= (x^{(0,1)})^2 + x^{(1,0)} \in \overline{\mathbb{F}_2}\{x\}.
\]
With analogous computations, for $\omega$ such that $(1,3) = \min_{\preceq_\omega} V$ we obtain
\[
    \initial_{(S, \omega)}(F)= x^{(2,3)} + x^{(1,0)},
\]
and for $\omega$ such that $(2,1) = \min_{\preceq_\omega} V$ we get 
\[
    \initial_{(S, \omega)}(F)= (x^{(0,1)})^2  + x^{(2,3)} + x^{(1,0)}.
\]
This reflects exactly the fact that the weighted vertices $(3,0)$ and $(1,3)$ in $\trop_v(F)(S)$ are attained twice at, respectively, the monomials $(x^{(0,1)})^2$ and $x^{(1,0)}$, and $x^{(2,3)}$ and $x^{(1,0)}$, while $(2,1)$ is attained three times, at all the monomials in the support of $F$.
\end{ex}

\subsection*{A differential version of \cite[Theorem 4.2]{draisma} and \cite[Proposition 2.2]{payneanal}}
    The statement of  \cite[Theorem 4.2]{draisma} and \cite[Proposition 2.2]{payneanal} is sometimes seen as a fourth part in the statement of the fundamental theorem of tropical geometry. Before moving to the next section we want to remark that a differential version of holds, and sketch a proof here. Let $(R^K_m,\{d/dt_i\})$ be equipped with the differentially enhanced valuation $\mathbf{v}$ to the pair $\mathbf{\Gamma}_m^K$ and let $I \subset \diff{R^K_m}{n}$ be a differential ideal. Denote as $A$ the finitely generated differential algebra $\diff{R^K_m}{n}/I$ over $R^K_m$ and let $\mathbf T$ be a $\mathbf{\Gamma}_m^K$-algebra. A differentially enhanced valuation $\mathbf{w}=(\widetilde{w}, w) \co A \rightarrow \mathbf T$ is said to be compatible with $\mathbf{v}$ if the diagram
	\begin{equation} \label{diagr:comp-diff-enh}
	\begin{tikzcd}
	&  &   (\widetilde{\Gamma}_K)_{m,v_K} \arrow[r] \arrow[dd] & T_1 \arrow[dd] \\ 
	R_m^K  \arrow[rrd, "v"'] \arrow[rru, "\widetilde{v}"] \arrow[r] &  A \arrow[rrd, crossing over, near end, "w"] \arrow[rru, crossing over, "\widetilde{w}"', near end]                         &               & \\
	& & \textup{Conv}_m(\widetilde{\Gamma}_K) \arrow[r]            & T_0           
	\end{tikzcd}
	\end{equation}
	commutes.
	Similarly to the non-differential case, the \emph{differential Berkovich space over $\mathbf{T}$} of $A$ is the set $\mathit{Berk}_{\mathbf{T}}(A)$ of differentially enhanced valuations $\mathbf{w}:=(w, \widetilde w) \co A \to \mathbf{T}$ that are compatible with $\mathbf{v}$.
	
	It is possible to define a differential tropicalization functor $\DTrop (-)$ (taking a differential presentation $\varphi$ of $A$ to an $\mathbf S$-algebra, see \cite[Section 5]{gianmereta})  and a universal differential presentation of $A$. With these definitions, a differential version of Payne's inverse limit theorem and of \cite[Theorem 4.4.1]{univtrop} holds, see \cite[Corollary 5.5.1]{gianmereta}. As morphisms $\DTrop (\varphi) \rightarrow  \mathbf T$ can be identified with a set of tuples in $T_1$ with specific properties (see \cite[Proposition 4.4.1]{gianmereta}), it is not hard to prove that 
		\begin{equation} \label{eq:4th-fund}
		\textup{Hom}_{\alg {\mathbf S}}(\DTrop(\varphi), \mathbf T) = \{(\widetilde w (\varphi(x_1)), \dots , \widetilde w (\varphi(x_n))) \in T_1^n \mid (w, \widetilde w) \in \mathit{Berk}_{\mathbf{T}}(A) \}.
		\end{equation}
	Moreover, when $\mathbf T = \mathbf \Gamma^K_m$, we obtain an equality between the right hand side of Equation \ref{eq:4th-fund} above and the left hand side of Theorem \ref{theorem:fundamental-initial}.

\section{The radius of convergence of solutions to a nonarchimedean differential equation}\label{section:raggi}
In this section we focus on the case $m=1$ and introduce a notion of radius of convergence for elements of $\pt$. Then, given a nontrivially valued field
$K$, we state a corollary of Theorem \ref{theorem:fundamental} relating this tropical notion of radius of convergence to the classical one for solutions of a differential equation over $\pk$.

\begin{defin}\label{definition:conv-radius}
    Given $A = \sum_{i=0}^\infty a_i t^i \in \pt$, we define its radius of convergence with respect to the real number $c > 1$ as:
    \[
        r_c(A):= \sup \{r \in [0,\infty) \mid \lim_{i \to +\infty} c^{-a_i} r^i = 0\} \in [0,\infty]
    \]
    with the convention that $c^{- \infty} = 0$ for any $c > 1$.
\end{defin}

\begin{oss}
    Notice that choosing any other real number $c'$ in the definition above gives the following relation for a tropical power series $A$:
    \[
        r_c(A) = r_{c'}(A) ^{\frac{1}{\log_c(c')}}.
    \]
\end{oss}
Let us consider now a valued field $(K,v_K)$, and assume the valuation $v_K$ has a section. For every $x \in K$, the norm on $K$ associated to $v_K$ is defined as $|x|_{K}:=c ^{-v_K(x)}$ for some $c > 1$. Canonical choices for $c$ are $c:=p$ when $v$ is the $p$-adic valuation for some prime number $p$ and $c:=e$ when $v$ is the $t$-adic valuation. 

Given a power series $A = \sum_{i=0}^\infty a_i t^i  \in \pk$ its radius of convergence is defined as follows: 
\[
        r(A):= \sup \{r \in [0,\infty) \mid \lim_{i \to +\infty} |a_i|_{K} \, r^i = 0\} \in [0,\infty].
\]
It is clear that, endowing $\pk$ with the differentially enhanced valuation $\mathbf{v}$, for every $A \in \pk$ we have the equality:
\[
    r(A) = r_{c}(\trop_{\widetilde v}(A)).
\]
where $c$ is chosen to be the same as in the definition of $| \! -  \! |_{K}$. 
Thus the following corollary of the Fundamental Theorem follows directly:

\begin{cor}\label{corolla:radius}
    Let $I \subset \pk\{x\}$ be a differential ideal and $r \in \{r_{c}(S)  \mid S \in \textup{Sol}_{\Gamma^K_m}(\trop_v(I))\}$. There exist a valued field extension $(L,v_L)$ of $(K,v_K)$ and $A \in \textup{Sol}_{L[\![t]\!]}(I)$ such that $r(A) = r$.
\end{cor}
In future work, we plan to determine hypotheses under which the field extension is not needed for Theorem \ref{theorem:fundamental-introduction} to hold, and to have a more explicit control over the field extension.

\subsection*{An explicit computation} To conclude, we give here an example of computation of the solutions of the tropicalization of a linear $p$-adic differential equation.
\begin{ex}\label{ex:p-adic-diffeq}
Let $p$ a prime number and consider the differential enhancement $\mathbf v$ as above for $K=\C_p$ endowed with the $p$-adic valuation. Let $\zeta \in \C_p$ such that $\zeta^{p-1} = -p$, and consider the equation $f = x' - p \zeta t^{p-1} x  \in \C_p[\![t]\!]\{x\}$. $p$-adic solutions at the origin to this equations are $p$-adic multiples of $\exp( \zeta t^p) = \sum_{n =0}^\infty (1/n!) (\zeta t^p)^n \in \C_p[\![t]\!]$. The interest of this $p$-adic differential equation lies in the fact that it is somehow the easiest example for which we can observe the piecewise linear behaviour of the radius of convergence as a function of the distance of the expansion point from the origin. In particular, $\exp( \zeta t^p)$ has radius of convergence 1 but there exists an $r > 1$ such that for solutions at $|x|_p > r$ the radius of convergence decreases. Let us see a concrete example of a computation of the solutions to a linear tropical differential equation, by computing the set of solutions to the tropicalization of the generators (as an ideal) of the differential ideal generated by $f$. 

One can prove by induction that $n$-th derivative of $f$ can be expressed as: 
\begin{equation*}
d^n f =
\begin{cases*}
x^{(n+1)} - \sum_{i=0}^n  \binom{n}{i} \zeta \prod_{k=0}^{n-i} (p-k) t^{p-1-n+i} x^{(i)} & if $n < p$ \\
x^{(n+1)} - \sum_{i=0}^{p-1} \binom{n}{p-1-i} \zeta \prod_{k=0}^{p-1-i} (p-k) t^{p-1-i} x^{(i+n-p+1)}  & if $n \ge p$
\end{cases*}
\end{equation*}
which implies: 
\begin{equation*}
\trop_v(d^n f ) =
\begin{cases*}
x^{(n+1)} + \sum_{i=0}^n (p-1-n+i, v_p( \binom{n}{i}) + \frac{p}{p-1})  x^{(i)}& if $n < p$ \\
x^{(n+1)} + \sum_{i=0}^{p-1}  (i, v_p( \binom{n}{p-1-i}) + \frac{p}{p-1})  x^{(i+n-p+1)}    & if $n \ge p$
\end{cases*}
\end{equation*}
We prove that a solution $S= \sum_{k=0}^\infty s_k t^k$ to the tropical system $\{\trop_v(d^n f)\}_{n \in N}$ is of the form $c \odot A$, where $c \in \T$ and $A = \sum_{k=0}^\infty a_k t^k \in \pt_{v_p}$ for:
\begin{equation*}
a_k =
\begin{cases*}
\infty & if $p \ndiv k$ \\
\frac{m}{p-1} - v_p(m!)  & if $k = mp$
\end{cases*}
\end{equation*}
i.e.\ that $A$ is the tropicalization of $\exp (\zeta t^p)$ and the set of solutions to $\trop_v(I)$ is a tropical linear space as expected. We prove this by induction on $k$. For $k=0$ this is trivially satisfied as $p \mid 0$ and $0 = \frac{0}{p-1} - v_p(0!)$. For $0 < k \le p$, consider $\textup{trop}_v(f) = x' + \left(p-1, \frac{p}{p-1}\right)x$ then we have that
\[
\textup{trop}_v(f)(A) = (0, a_1) \oplus \left (p-1, \frac{p}{p-1} \right) ,
\]
and as we want $A$ to be a solution, i.e.\ the minimum to be attained twice in the equation above, $a_k$ has to be equal to $\infty$ for every $0 < k < p$, thus we obtain equality in the first coordinate: 
\[
\textup{trop}_v(f)(A) = (p-1, a_p +1) \oplus \left ( p-1, \frac{p}{p-1} \right ).
\]
Again, as the minimum has to be attained twice in order for $A$ to be a solution, we obtain $a_p = \frac{p}{p-1}-1 =\frac{1}{p-1} = \frac{1}{p-1} - v_p(1!) $.

In general, let  $Mp < k \le (M+1)p$ for some $M \in \N$ and assume the inductive hypothesis holds for $Mp$. Notice that $v_p \left ( \binom{Mp}{p-1-i} \right )= 1$ for all $M \in \N$ and for all $i \in \{0, \dots , p-2\}$, thus:
\[
\trop_v(d^{Mp}f) = x^{(Mp+1)} + \left ( p-1, \frac{p}{p-1} \right ) x^{(Mp)} +  \sum_{j=1}^{p-1}  \left ( p-1-j, \frac{p}{p-1} +1 \right )  x^{(Mp- j)}
\] 
where $j:=i-p+1$.  We obtain that $\trop_v(d^{Mp}f) (A)$ is equal to:
\begin{equation*}
\begin{split}
& \left ( 0, a_{Mp+1} + \sum_{m=1}^M v_p(mp) \right ) \oplus   \left ( p-1, \frac{p}{p-1} \right ) \odot \left  ( 0, \frac{M}{p-1} - v_p(M!) +  \sum_{m=1}^M v_p(mp) \right )  \oplus \\
&  \bigoplus_{j=1}^{p-1}  \left ( p-1-j, \frac{p}{p-1} +1 \right ) \odot \left (j  ,   \frac{M}{p-1} - v_p(M!) + \sum_{m=1}^{M} v_p(mp)  \right )\\
\end{split}
\end{equation*}
As the first coordinate of every summand of the sum above, except for the first one, is equal to $p-1$ and $A$ is a solution, we have as before that $a_k$ has to be equal to $\infty$ for $Mp < k < (M+1)p$. For $k=(M+1)p$, since the following holds:
\[
\sum_{m=1}^M v_p(mp) = \sum_{m=1}^M ( 1+ v_p(m) ) = M + v_p(M!)
\]
we have:
\begin{equation*}
\begin{split}
\trop_v(d^{Mp}f) (A) = & \left ( p-1, a_{(M+1)p} + M+1 +v_p((M+1)!) \right ) \oplus    \left  ( p-1 , \frac{M +p}{p-1} + M  \right )  \oplus \\
&  \oplus \left ( p-1 , \frac{M +p}{p-1}   + M +1 \right )\\
\end{split}
\end{equation*}
which gives:
\[
a_{(M+1)p} +1 +v_p((M+1)!) = \frac{M +p}{p-1} .
\]
The last equality implies $a_{(M+1)p}= \frac{M +1}{p-1} - v_p((M+1)!)$. Lastly, it is easy to check that every tropical multiple of $A$ is again a solution to the tropical system $\{\trop_v(d^n f)\}_{n \in N}$, as these are linear equations.

\begin{oss}
It follows from the computation above that no field extension is needed in this case, even if the differential polynomial $f$ does not satisfy the hypothesis of part (2) of Proposition \ref{Prop:limit-empty-iff-some-empty}.
\end{oss}

The radius of convergence of a tropical solution $A\neq \infty$ to $f$ is equal to
\[
    r_p(A) = \sup \{r \in [0,\infty) \mid \lim_{ m \rightarrow \infty} p^{-\frac{m}{p-1}} p^{v_p(m!)} r^{mp} = 0\}
\]
From Legendre's formula for $v_p(m!)$, it follows that the radius of convergence of the $p$-adic exponential is $p^{-\frac{1}{p-1}}$, thus $r_p(A) = 1$. From Corollary \ref{corolla:radius} this is true also for every nonzero solution of $f$.
\end{ex}
    \begin{oss}
    Notice that in Example \ref{ex:p-adic-diffeq} above we also proved that $\{f\} \subset I$ is a \emph{tropical differential basis} for $I$, in the sense of \cite[Definition 4.1]{FT20}, i.e.\ a set $\mathcal G \subset I$ such that 
    \[
        \textup{Sol}_{\mathbf{S}}(\trop_v(I)) = \bigcap_{g \in \mathcal G} \bigcap_{k \in \N}\textup{Sol}_{\mathbf{S}}(\trop_v(d^k g)).
    \]
    In the same paper the authors, working in Grigoriev's setting, give an example of a differential ideal $I$ generated by linear forms that does not admit a finite tropical differential basis of linear forms, in contrast with the behaviour of tropical basis in the non-differential setting (see \cite[Theorem 2.6.6]{macsturm}).
    \end{oss}
    \subsection*{Conclusions and future work} The fundamental theorem hereby proved in full generality opens several directions of research. This work takes the first steps towards the development of methods using tropical techniques to compute the radius of convergence function, without a priori knowing the solutions to the differential equation. 
    
    To this aim, we believe the following research directions should be undertaken: 
    \begin{itemize} 
    \item Characterizing the class of differential ideals that require a field extension of the coefficients and control this extension, in particular its value group;
    \item Introducing a notion of tropical differential basis satisfying a finiteness condition as in the classical case, or at least sufficient conditions on the equations for finiteness of basis as above to hold. In the linear case this would make similar computation as in Example \ref{ex:p-adic-diffeq} possible, without prior knowledge on the solutions over $\pk$;
    \item The methods introduced in this work are only effective for computing single values of the radius of convergence function. In order to capture its global behaviour, it is necessary to develop techniques to move the tropical point of expansion. This could be done via a theory of tropical differential modules coherent with the $p$-adic for tropicalized equations;
    \end{itemize}
    Lastly, nothing is known about the polyhedral geometry of the set of solutions to a system of TDEs, not even when it comes from a classical system. To this end, the third item in Theorem \ref{theorem:fundamental-introduction} is foundational.

\section*{Table of notations}
If not specified otherwise, the following notations will have the following meaning:
\begin{center}
\begin{longtable}{p{3cm} p{12cm}}
$m$ & a positive integer \\
$[m]$ & the set $\{1, \dots , m\}$\\
$e_i$ & the $i$-th standard basis vector of $\R^m$ \\
$\langle \cdot, \cdot \rangle $ & the scalar product of two vectors of $\R^m$ \\
$\tau, \rho$ & a strongly convex, rational polyhedral cone of dimension $m$ in $\R^m$\\
$C_\tau(Z)$ &the convex hull with respect to $\tau$ of a subset $\Z \subseteq \R^m$  \\
$\textup{Vert}_\tau(Z)$ & the set of vertices of $C_\tau(Z)$\\
$\B$ & the idempotent semiring $(\{0, \infty\}, \min, +)$\\
$\T$ & the idempotent semiring $(\R\cup \{\infty\}, \min, +)$\\
$\T^{(m)}$ & the idempotent semiring $(\R^m\cup \{\infty\}, \min_{\text{lex}}, +)$ \\
$\widetilde{G}$ & the idempotent semiring $(G \cup \{\infty\}, \min_\preceq, +)$ for a totally ordered monoid $(G,+, \preceq)$\\
$\omega$ & a vector in $(\R_{>0})^m$ with rationally independent coordinates\\
$H_\omega$ & the halfspace $\{x \in \R^m \mid \langle x, \omega \rangle \ge 0 \}$\\
$\textup{Lead}_\omega(G, S)$ & the idempotent semiring of leading terms with respect to $\omega$ of an additive submonoid $G \subseteq \Q^m$, weighted in $S$\\
$\textup{Conv}(\tau,S)$ &  the idempotent semiring of convex hulls w.r.t.\ $\tau$ with vertices weighted in an idempotent semiring $S$ \\
$\textup{Conv}_m(S)$ &  the idempotent semiring $\textup{Conv}(\tau,S)$ for $\tau = (\R_{\ge 0})^m$ \\
$\textup{Conv}(\tau,S)_{\{e_i\}}$ & the localization of $\textup{Conv}(\tau,S)$ in the translates of $\tau$\\
$\sigma_0$ & the support map $\textup{Conv}(\tau,S) \rightarrow \textup{Conv}(\tau,\B)$ of Remark \ref{oss:Conv-m=1} \\
$K$ & an algebraically closed field of characteristic 0 \\
$v_\text{triv}$ & the trivial valuation $K \rightarrow \B$\\
$v_K$ & a nontrivial valuation $K \rightarrow \T$\\
$R_K$ & the valuation ring of $v_K$ \\
$\mathfrak m_K$ & the maximal ideal of $R_K$ \\
$k$ & the residue field $R_K / \mathfrak m _K$ of $v_K$ \\
$\pi_K$ & a uniformizer of $v_K$ \\
$\Gamma_K$ & the value group of $v_K$ \\
$R_\nu$ & the valuation ring of a valuation $\nu$ \\
$\mathfrak m_\nu$ & the maximal ideal of $R_\nu$, for a valuation $\nu$\\
$\supp(A)$ & the support of the power series A\\
$\pkcone \tau$ & the ring of power series over $K$ with support in $ \Z^m \cap \tau$\\
$R_m^K$ & the ring $\pkk m$, i.e.\ the ring $\pkcone \tau$ for $\tau = (\R_{\ge 0})^m$\\
$\pkcone \tau _{t_1, \dots, t_m}$ & the localization of $\pkcone \tau$ in the variables\\
$K(\!(t_1, \dots , t_m)\!)_\omega$ & the field of multivariate Laurent series with support contained in a translate of $H_\omega$ \\
$w_\tau$ & the valuation $\pkcone \tau \rightarrow \textup{Conv}(\tau, \B)$ of point (1) of Example \ref{ex:conv-vals} \\
$v_\tau$ & the valuation $\pkcone \tau \rightarrow \textup{Conv}(\tau, \T)$ of point (2) of Example \ref{ex:conv-vals} \\
$v_\omega$ & the valuation $K(\!(t_1, \dots , t_m)\!)_\omega \rightarrow \textup{Lead}_\omega(\Z^m, \T)$ as in \ref{eq:v-omega} \\
$w$ & the valuation $w_\tau$ for $\tau = (\R_{\ge0})^m$ \\
$v$ & the valuation $v_\tau$ for $\tau = (\R_{\ge0})^m$ \\
$\B_1$ & the strict differential semiring $\pb$ of Boolean power series \\
$\T_{1,v}$ &  the (non-strict) differential semiring $\pt$ of tropical power series, with differential w.r.t.\ a nontrivial valuation $v \co \N \rightarrow \T$  \\
$\B_m$ &  the strict partial differential semiring $\pbb m$ of multivariate Boolean power series \\
$\T_{m,v}$ &  the (non-strict) partial differential semiring $\ptt n _v$ of multivariate tropical power series, with differentials w.r.t.\ a nontrivial valuation $v \co \N \rightarrow \T$\\
$\diff R n$ & the differential algebra of differential polynomials in $n$ variables over a differential ring $R$ as in \cite{ritt} \\
$\diff S n$ & the differential algebra of differential polynomials in $n$ variables over a differential semiring $S$ as in \cite{gianmereta} \\
$\basic S n$ & the $S$-algebra $S[x_i^{(\jj)} \mid i \in [n], \jj \in \N^m]$ over a differential semiring $S$ with $m$ differentials \\
$\diff{\mathbf S} n$ & the $\mathbf S$-algebra of differential polynomials in $n$ variables over a (partial) tropical pair $\mathbf S$, as in \cite{gianmereta}\\
$\mathbf T _m$ & the (partial) tropical pair $\Psi_m \co \B_m \to \textup{Conv}_m(\B)$ of point (3) of Example \ref{ex:pairs} \\
$\mathbf T$ & the tropical pair $\mathbf T _m$ for $m=1$, i.e.\  $\Psi \co \B_1  \to \T$ \\
$\mathbf S _m$ & the (partial) tropical pair $\Phi_m \co \T_{m,v} \to \textup{Conv}_m(\T)$ of point (4) of Example \ref{ex:pairs} \\
$\mathbf S$ & the tropical pair $\mathbf S _m$ for $m=1$, i.e.\ $\Phi \co \T_{1,v}  \to \T^{(2)}$  \\
$\mathbf \Gamma^K_{m}$ & the (partial) tropical pair obtained from $\mathbf S _m$ restricting coefficients and weights to $\Gamma_K$, as in Remark \ref{oss:restr-of-diff-enh} \\
$\sigma$ & the morphism of pairs $(\sigma_0, \sigma_1 ) \co \mathbf S _m \rightarrow \mathbf T _m$, as in Diagram \ref{diagram:refined-to-grig}\\
$\widetilde w$ & submultiplicative seminorm $R_m \rightarrow \B_m$ given by coefficientwise trivial valuation, as in point (3) of Example \ref{ex:diff-enhancements} \\
$\widetilde v$ & submultiplicative seminorm $R_m \rightarrow \T_{m,v}$ given by coefficientwise application of a valuation $v$, as in point (4) of Example \ref{ex:diff-enhancements} \\
$\mathbf w$ & the differentially enhanced valuation $(w, \widetilde w) \co R_m \rightarrow \mathbf T _m$ of point (3) of Example \ref{ex:diff-enhancements}, used in \cite{grig, aroca, sebastian, FT20} \\
$\mathbf v$ & the differentially enhanced valuation $(v, \widetilde v)  \co R_m \rightarrow \mathbf S _m$ of point (4) of Example \ref{ex:diff-enhancements}, introduced in \cite{gianmereta, tesimereta} 
\end{longtable}
\end{center}

\newcommand{\etalchar}[1]{$^{#1}$}
	
\end{document}